\documentclass[11pt,reqno]{amsart}
\usepackage{amsmath,amsopn,amssymb,amsthm,multicol}
\usepackage{hyperref}
\usepackage{graphicx}

\DeclareMathOperator{\Ad}{Ad}

\DeclareMathOperator{\Ric}{Ric}

\DeclareMathOperator{\Mark}{Mrk}

\newcommand{\fr}{\mathfrak}
\newcommand{\al}{\alpha}
\newcommand{\be}{\beta}

\newcommand{\bb}{\mathbb}

\DeclareMathOperator{\SO}{SO}
\DeclareMathOperator{\s}{S}

\DeclareMathOperator{\Sp}{Sp}
\DeclareMathOperator{\SU}{SU}
\DeclareMathOperator{\U}{U}
\DeclareMathOperator{\G}{G}

\DeclareMathOperator{\E}{E}

 \newtheorem{lemma} {Lemma} [section]
\newtheorem{theorem}[lemma]{Theorem}

\begin{document}

\title{Progress on homogeneous Einstein manifolds and some open probrems} 
\author{Andreas Arvanitoyeorgos}
\address{University of Patras, Department of Mathematics, GR-26500 Patras, Greece}
\email{arvanito@math.upatras.gr}

   \begin{abstract}
We give an overview of progress on homogeneous Einstein metrics on large classes of homogeneous manifolds, such as generalized flag manifolds and  Stiefel manifolds. The main difference between these two classes of homogeneous spaces is that their isotropy representation does not contain/contain equivalent summands.  We also discuss a third class of homogeneous spaces that falls into the intersection of such dichotomy, namely the generalized Wallach spaces.  We give new invariant Einstein metrics on the Stiefel manifold $V_5\mathbb{R}^n$ ($n\ge 7$) and through this example we show how to prove existence of invariant Einstein metrics by manipulating parametric systems of polynomial equations.  This is done by using Gr\"obner bases techniques.
Finally, we discuss some open problems.

  \medskip
\noindent 2010 {\it Mathematics Subject Classification.} Primary 53C25; Secondary 53C30, 22C05, 22E60, 13P05, 13P10, 13P15

\medskip
\noindent {\it Keywords}:    Homogeneous space; Einstein metric; isotropy representation; generalized flag manifold; Stiefel manifold;  generalized Wallach space; algebraic system of equations; 
Gr\"obner basis
   \end{abstract}

\maketitle


 \section{Introduction}
\markboth{Andreas Arvanitoyeorgos}{Recent progress on homogeneous Einstein manifolds and some open problems}

A Riemannian manifold $(M, g)$ is called Einstein if it has constant Ricci curvature, i.e. $\Ric_{g}=\lambda\cdot g$ for some $\lambda\in\bb{R}$. 
Due to an old result of Hilbert, for $M$ compact, Einstein metrics are precisely the critical points of the scalar curvature functional over the set of Riemannian metrics of volume $1$.
If $\lambda > 0$ then $M = G/H$ is compact,  if
 $\lambda = 0$, $M$ is Ricci flat, and if 
$\lambda < 0$ then $G/H$ non-compact.
  For results on Einstein manifolds before 1987 we refer to the book by A. Besse \cite{Be}. 
  The two  articles \cite{W1}, \cite{W2} of M. Wang  contain results up to 1999 and 2013 respectively.
  General existence results are difficult to obtain and some methods are described
  in \cite{Bom}, \cite{BWZ} and \cite{WZ}.  
  
For homogeneous spaces $M=G/K$ the problem is to find and classify all $G$-invariant Einstein metrics,
or decide if the set of $G$-invariant Einstein metrics is finite or not.
A conjecture by W. Ziller (\cite{BWZ}) states that  if ${\rm rk}(G) = {\rm rk}(H)$, or 
if the isotropy representation of $G/K$ contains only non equivalent summands,
then the number of $G$-invariant Einstein metrics is finite.  

The aim of the present article is to give an overview of recent progress on homogeneous Einstein metrics on large classes of homogeneous spaces, namely generalized flag manifolds, Stiefel manifolds and generalized Wallach spaces. 
The first two  classes of homogeneous spaces fall into a certain dichotomy between homogeneous spaces whose isotropy representation decomposes into non equivalent summands and those whose isotropy representation contain equivalent summands.  Generalized Wallach spaces lie in both classes of spaces.  Furthermore, interesting and non trivial problems arise in both classes of spaces, which are discussed at the end of the paper as open problems.

The article also contains some  new results concerning existence of new invariant Einstein metrics on the Stiefel manifold $V_5\mathbb{R}^n$ $(n\ge 7)$, thus extending a result in \cite{ASS}.  
Through this class of homogeneous spaces we have the opportunity to show how to prove existence on invariant Einstein metrics by manipulating parametric systems of polynomial equations.  This is done by using Gr\"obner bases techniques.

\section{The Ricci tensor for reductive homogeneous spaces}
In this section we recall an expression for the Ricci tensor for an $G$-invariant Riemannian
metric on a reductive homogeneous space whose isotropy representation
is decomposed into a sum of non equivalent irreducible summands.

Let $G$ be a compact semisimple Lie group, $K$ a connected closed subgroup of $G$  and  
let  $\frak g$ and $\fr{k}$  be  the corresponding Lie algebras. 
The Killing form $B$ of $\frak g$ is negative definite, so we can define an $\mbox{Ad}(G)$-invariant inner product $-B$ on 
  $\frak g$. 
Let $\frak g$ = $\frak k \oplus
\frak m$ be a reductive decomposition of $\frak g$ with respect to $-B$ so that $\left[\,\frak k,\, \frak m\,\right] \subset \frak m$ and
$\frak m\cong T_o(G/K)$.
 We assume that $ {\frak m} $ admits a decomposition into mutually non equivalent irreducible $\mbox{Ad}(K)$-modules as
\begin{equation}\label{iso}
{\frak m} = {\frak m}_1 \oplus \cdots \oplus {\frak m}_q.
\end{equation} 
Then any $G$-invariant metric on $G/K$ can be expressed as  
\begin{eqnarray}
 \langle\  , \  \rangle  =  
x_1   (-B)|_{\mbox{\footnotesize$ \frak m$}_1} + \cdots + 
 x_q   (-B)|_{\mbox{\footnotesize$ \frak m$}_q},  \label{eq2}
\end{eqnarray}
for positive real numbers $(x_1, \dots, x_q)\in\bb{R}^{q}_{+}$.  Note that  $G$-invariant symmetric covariant 2-tensors on $G/K$ are 
of the same form as the Riemannian metrics (although they  are not necessarilly  positive definite).  
 In particular, the Ricci tensor $r$ of a $G$-invariant Riemannian metric on $G/K$ is of the same form as (\ref{eq2}), that is 
 \[
 r=y_1 (-B)|_{\mbox{\footnotesize$ \frak m$}_1}  + \cdots + y_{q} (-B)|_{\mbox{\footnotesize$ \frak m$}_q} ,
 \]
 for some real numbers $y_1, \ldots, y_q$.

Let $\lbrace e_{\alpha} \rbrace$ be a $(-B)$-orthonormal basis 
adapted to the decomposition of $\frak m$,    i.e. 
$e_{\alpha} \in {\frak m}_i$ for some $i$, and
$\alpha < \beta$ if $i<j$. 
We put ${A^\gamma_{\alpha
\beta}}= -B \left(\left[e_{\alpha},e_{\beta}\right],e_{\gamma}\right)$ so that
$\left[e_{\alpha},e_{\beta}\right]
= \displaystyle{\sum_{\gamma}
A^\gamma_{\alpha \beta} e_{\gamma}}$ and set 
$\displaystyle{k \brack {ij}}=\sum (A^\gamma_{\alpha \beta})^2$, where the sum is
taken over all indices $\alpha, \beta, \gamma$ with $e_\alpha \in
{\frak m}_i,\ e_\beta \in {\frak m}_j,\ e_\gamma \in {\frak m}_k$ (cf. \cite{WZ}).  
Then the positive numbers $\displaystyle{k \brack {ij}}$ are independent of the 
$B$-orthonormal bases chosen for ${\frak m}_i, {\frak m}_j, {\frak m}_k$,
and 
$\displaystyle{k \brack {ij}}\ =\ \displaystyle{k \brack {ji}}\ =\ \displaystyle{j \brack {ki}}.  
 \label{eq3}
$

Let $ d_k= \dim{\frak m}_{k}$. Then we have the following:

\begin{lemma}\label{ric2}\textnormal{(\cite{PS})}
The components ${ r}_{1}, \dots, {r}_{q}$ 
of the Ricci tensor ${r}$ of the metric $ \langle  \,\,\, , \,\,\, \rangle $ of the
form {\em (\ref{eq2})} on $G/K$ are given by 
\begin{equation}
{r}_k = \frac{1}{2x_k}+\frac{1}{4d_k}\sum_{j,i}
\frac{x_k}{x_j x_i} {k \brack {ji}}
-\frac{1}{2d_k}\sum_{j,i}\frac{x_j}{x_k x_i} {j \brack {ki}}
 \quad (k= 1,\ \dots, q),    \label{eq51}
\end{equation}
where the sum is taken over $i, j =1,\dots, q$.
\end{lemma} 
Since by assumption the submodules $\fr{m}_{i}, \fr{m}_{j}$ in the decomposition (\ref{iso}) are mutually non equivalent for any $i\neq j$, it is $r(\fr{m}_{i}, \fr{m}_{j})=0$. 
If $\fr{m}_i\cong\fr{m}_j$ for some $i\ne j$ then we need to check whether 
$r(\fr{m}_{i}, \fr{m}_{j})=0$. This is not an easy task in general.
Once the condition $r(\fr{m}_{i}, \fr{m}_{j})=0$ is confirmed we can use
  Lemma \ref{ric2}.   Then    $G$-invariant Einstein metrics on $M=G/K$ are exactly the positive real solutions $g=(x_1, \ldots, x_q)\in\bb{R}^{q}_{+}$  of the  polynomial system $\{r_1=\lambda, \ r_2=\lambda, \ \ldots, \ r_{q}=\lambda\}$, where $\lambda\in \bb{R}_{+}$ is the Einstein constant.
  If some of the submodules in (\ref{iso}) are equivalent as $\Ad(K)$-modules, then the computation of the Ricci tensor is more laborious (e.g. using basic formulas from \cite{Be} or proving other variations of these).

 \section{Homogeneous spaces with non equivalent isotropy summands}

A generalized flag manifold is a a homogeneous space $M=G/K$ where $G$ is a compact semisimple Lie group and $K$ is the cetralizer of a torus in $G$.  Equivalently, it is diffeomorphic to the adjoint orbit $\Ad(G)w$, for some $w\in\frak g$, the Lie algebra of $G$. Typical examples are the manifolds of partial flags  $\SU(n)/\s(U(n_1)\times\cdots\times \U(n_k))$ ($n=n_1+\cdots +n_k$) and full flags $\SU(n)/T$ in $\mathbb{C}^n$, where $T=\s(\U(1)\times\cdots\U(1))$ is a maximal torus in $\SU(n)$.

Next, we will present the Lie theoretic description of a generalized flag manifold.
\subsection{Decomposition associated to generalized flag manifolds}
Let $G$ be a compact semisimple Lie group,  
$\frak g$ the Lie
algebra of $G$ and $\frak h$ a maximal abelian subalgebra of
$\frak g$. We denote by ${\frak g }^{\mathbb C}_{}$  and ${\frak h}^{\mathbb C}_{}$ the complexification of $\frak g$ and 
$\frak h$
respectively. 
  We identify an element of the root system $\Delta$ of
 ${\frak g }^{\mathbb C}_{}$ relative to the Cartan subalgebra  
 ${\frak h}^{\mathbb C}_{}$ with an element of ${\frak h}_0 =  \sqrt{-1}\frak h$ by the
duality defined by the Killing form of ${\frak
g}^{\mathbb C}_{}$. Let $\Pi$ = $\{\alpha^{}_1, \dots, \alpha^{}_l\}$
be a fundamental system of $\Delta$ and $\{\Lambda^{}_1, \dots,
\Lambda^{}_l\}$ the fundamental weights of ${\frak g }^{\mathbb C}_{}$ 
corresponding to $\Pi$, that is 
$\frac{2(\Lambda^{}_i, \alpha^{}_j)}{(\alpha^{}_j, \alpha^{}_j)} =
\delta^{}_{ij},\ (1 \le i, j \le l).$ 
Let $\Pi^{}_0$ be a subset of $\Pi$ and $\Pi - \Pi^{}_0$ =
$\{\alpha^{}_{i_1}, \cdots, \alpha^{}_{i_r}\}$, where  $1 \le
{i_1} < \cdots <{i_r} \le l $. We put
$[\Pi^{}_0] = \Delta\cap\{\Pi^{}_0\}^{}_{\mathbb Z}$, where
 $\{\Pi^{}_0\}^{}_{\mathbb Z}$ denotes the subspace of ${\frak h}_0$ generated by $\Pi^{}_0$. 
Consider the root space decomposition of  ${\frak g }^{\mathbb C}_{}$
relative to  ${\frak h}^{\mathbb C}_{}$ as
$ {\frak g }^{\mathbb C}_{}  =  {\frak h}^{\mathbb C}_{} +
\sum^{}_{\alpha \in \Delta} {\frak g }^{\mathbb C}_{\alpha} .$
We define a parabolic subalgebra  ${\frak u}^{}_{}$ of  
${\frak g }^{\mathbb C}_{}$ by 
$\displaystyle  {\frak u}^{}_{}  =  {\frak h}^{\mathbb C}_{} +
\sum^{}_{\alpha \in [\Pi^{}_0]\cup\Delta^+_{}}
{\frak g }^{\mathbb C}_{\alpha},$
where $\Delta^+_{}$ is the set of all positive roots relative to $\Pi$. 
Note that the nilradical ${\frak n}$ of $ {\frak u}$ is given by 
${\frak n}  = \sum^{}_{\alpha \in \Delta^+ - [\Pi^{}_0]} {\frak g }^{\mathbb C}_{\alpha} .
$ 
We put $ \Delta_{\frak m}^{+} =  \Delta^+ - [\Pi^{}_0]$. 

 Let $G^{\mathbb C}$ be a simply connected complex semi-simple
Lie group whose Lie algebra is ${\frak g }^{\mathbb C}_{} $ and $U$
the parabolic subgroup of $G^{\mathbb C}$ generated by ${\frak
u}^{}_{}$. Then the complex homogeneous manifold $G^{\mathbb C}/U$
is  compact simply connected and $G$ acts transitively on 
$G^{\mathbb C}/U$. Note also that $K = G\cap U$ is a connected closed subgroup
of $G$,  $G^{\mathbb C}/U$ = $G/K$ as $C^\infty_{}$-manifolds, and
$G^{\mathbb C}/U$ admits a $G$-invariant K\"ahler metric. 
Let $\frak k$ be the Lie algebra of $K$ and ${\frak k}^{\mathbb C}$ the complexification of $\frak k$. Then we have direct decompositions 
$\displaystyle {\frak u}^{}_{}  =  {\frak k}^{\mathbb C}_{} \oplus {\frak n}, 
\quad 
\displaystyle {\frak k}^{\mathbb C}_{}  =  {\frak h}^{\mathbb C}_{} +
\sum^{}_{\alpha \in [\Pi^{}_0]} {\frak g }^{\mathbb C}_{\alpha}.$
Take a Weyl basis $E_{-\alpha} \in {\frak g}^{\mathbb C}_{\alpha}
\,\,(\alpha \in \Delta )$ with
$$\begin{array}{ll} 
\left[E_{\alpha}, E_{-\alpha}\right] &=  -\alpha \,(\alpha \in \Delta)
 \\
{} \left[E_{\alpha}, E_{\beta}\right] &=  
\left\{
 \begin{array}{ll} N_{\alpha, \, \beta}E_{\alpha + \beta}& \mbox{if \quad} \alpha +\beta \in
   \Delta
   \\
   0 & \mbox{if \quad} \alpha +\beta \not\in \Delta, 
 \end{array}
\right.
\end{array}$$
where $N_{\alpha, \, \beta}$ = 
$N_{-\alpha, \, -\beta} \in {\mathbb R}.$
Then we have 
$${\frak g} = {\frak h} + \sum_{\alpha \in \Delta} \left\{
{\mathbb R}(E_{\alpha} + E_{-\alpha}) + {\mathbb R} \sqrt{-1} (E_{\alpha}
- E_{-\alpha})\right\}$$   
and the Lie subalgebra $\frak k$ is given
by 
$${\frak k} ={\frak h} + \sum_{\alpha \in \left[\Pi_0\right]}
 \left\{
{\mathbb R} (E_{\alpha} + E_{-\alpha}) + {\mathbb R}\sqrt{-1} (E_{\alpha}
- E_{-\alpha})\right\}.
$$ 
     
       For integers $j_1, \cdots, j_r$ with $(j_1, \cdots, j_r) \neq (0, \cdots, 0)$ , we put 
       $$\Delta( j_1, \cdots, j_r ) = \left\{ \  \sum_{j=1}^{l} m_j \alpha_j \in \Delta^{+} \ \  \Big\vert  \ \  m_{i_1} = j_1, \cdots, m_{i_r} = j_r  \ \right\}. 
       $$
       Note that $\displaystyle  \Delta_{\frak m}^{+} =  \Delta^+ - [\Pi^{}_0]= \bigcup_{j_1, \cdots, \, j_r} \Delta( j_1, \cdots, j_r )$.    
        
        For $  \Delta( j_1, \cdots, j_r ) \neq \emptyset$,  we define an $\mbox{Ad}_G(K)$-invariant subspace ${\frak m}( j_1, \cdots, j_r )$ of $\frak g$ by 
 $${\frak m}( j_1, \cdots,  j_r ) = \sum_{\alpha \in  \Delta( j_1, \cdots, \, j_r )}
 \left\{ {\mathbb R}(E_{\alpha} + E_{-\alpha}) + {\mathbb R} \sqrt{-1} (E_{\alpha} - E_{-\alpha})\right\}. $$
Then we have a decomposition of $\frak m$ into mutually non-equivalent irreducible $\mbox{Ad}_G(K)$-modules  $ {\frak m}( j_1, \cdots,  j_r )$ as
  $\frak m = \sum_{j_1, \cdots, \, j_r} {\frak m}( j_1, \cdots,  j_r ). 
  $ 
      We put ${\frak t} = \Big\{ H \in   {\frak h}_0 \  \vert \  ( H, \ \Pi_{0}) =   0   \Big\}$. 
      Then $\{ \Lambda^{}_{i_1}, \cdots, \Lambda^{}_{i_r}\}$ is a basis of ${\frak t}$. Put $ \frak s = \sqrt{-1}{\frak t} $. Then  the Lie algebra ${\frak k} $ is given by 
     ${\frak k} = {\frak z}({\frak s})$ (the Lie algebra of centralizer of a torus $S$ in $G$).  

 We consider the restriction map 
  $ \kappa  :  {\frak h}_0^* \to {\frak t}^*, \quad   \alpha \mapsto \alpha\vert_{\frak t}$
  and set  $\Delta_T = \kappa(\Delta)$.  The elements of $\Delta_T$ are called $T$-{\it roots}. 
  
   There exists (\cite{AlPe}) a 1-1 correspondence between $T$-roots $\xi$ and irreducible submodules $\displaystyle {\frak m}_{\xi}$ of the $\mbox{Ad}_G(K)$-module ${\frak m}^{\mathbb C}_{}$ given by 
 $$\displaystyle \Delta_T  \ni \xi \mapsto   {\frak m}_{\xi} = \sum_{\kappa(\alpha) = \xi} {\frak g }^{\mathbb C}_{\alpha}.$$ 
 Thus we have a decomposition of the $\mbox{Ad}_G(K)$-module ${\frak m}^{\mathbb C}_{}$ as 
   $ {\frak m}^{\mathbb C}_{} = \sum_{\xi \in \Delta_T} {\frak m}_{\xi}. $
 Denote by $\Delta_T^{+}$   the set of all positive $T$-roots, that is, the restricton of the system $\Delta^+$. Then 
we have $\displaystyle {\frak n}  = \sum_{\xi \in \Delta_T^{+}} {\frak m}_{\xi}$. 
Denote by $\tau$   the complex conjugation of ${\frak g }^{\mathbb C}_{}$ with respect to $\frak g$ (note that $\tau$ interchanges ${\frak g }^{\mathbb C}_{\alpha}$ and ${\frak g }^{\mathbb C}_{-\alpha}$) and by 
 ${\frak v}^{\tau}_{}$  the set of fixed points of $\tau$ in a complex vector subspace ${\frak v}$ of  ${\frak g }^{\mathbb C}_{}$.   Thus we have a decomposition of $\mbox{Ad}_G(K)$-module ${\frak m}$ into irreducible submodules as 
 ${ \frak m} = \sum_{\xi \in \Delta_T^{+}} \left( {\frak m}_{\xi} + {\frak m}_{-\xi} \right)^{\tau}.$ 
 
 There exists a natural 1-1 correspondence between $\Delta_T^{+}$ and the set $\left\{ \Delta( j_1, \cdots, j_r ) \neq \emptyset \right\}$. 
   For a generalized flag manifold $G/K$, we have a decomposition of $\frak m$ into  mutually non-equivalent   irreducible $\mbox{Ad}_G(H)$-modules  as
  $\displaystyle \frak m =  \sum_{\xi \in \Delta_T^{+}} \left( {\frak m}_{\xi} + {\frak m}_{-\xi} \right)^{\tau} =  \sum_{j_1, \cdots, \, j_r} {\frak m}( j_1, \cdots,  j_r ). 
  $ 
  Thus a $G$-invariant metric $g$ on $G/K$ can be written as  
\begin{eqnarray}
 g  =\sum_{\xi \in \Delta_T^{+}} x_{\xi} B|_{\left( {\frak m}_{\xi} + {\frak m}_{-\xi} \right)^{\tau}}  =  \sum_{j_1, \cdots, j_r}
x_{j_1 \cdots j_r}  B|_{ {\frak m}( j_1, \cdots,  j_r )}  \label{eq22}
\end{eqnarray}
for positive real numbers $ x_{\xi}$, $ x_{j_1 \cdots j_r}$.

   Put $ \displaystyle Z_{\frak t} = \Bigg\{\Lambda \in {\frak t} \ \Big\vert \ \frac{2(\Lambda, \ \alpha) }{( \alpha, \  \alpha)}  \in { \mathbb Z} \  \mbox{ for }  \mbox{ each } \alpha \in \Delta \Bigg\}$. 
Then  $ Z_{\frak t} $   is a lattice of  ${\frak t} $ generated by  $\{ \Lambda^{}_{i_1}, \dots, \Lambda^{}_{i_r}\}$.  
   For each $\Lambda \in  Z_{\frak t} $ there exists a unique holomorphic character $\chi_{\Lambda}$ of $U$
such that $\chi_{\Lambda}(\exp H) = \exp \Lambda( H ) $ for each $H \in  {\frak h}^{\mathbb C}_{}$. Then the correspondence $\Lambda \to \chi_{\Lambda} $
 gives an isomorphism of $Z_{\frak t}$ to the group of holomorphic characters
of $U$. 

Let $F_\Lambda$ denote the holomorphic line bundle on $G^{\mathbb C}/U$ associated to the
principal bundle $U \to G^{\mathbb C} \to  G^{\mathbb C}/U$ by the  holomorphic character $\chi_{\Lambda}$, and   $H( G^{\mathbb C}/U, \underline{\mathbb C}^*)$ the group of isomorphism classes of holomorphic line bundles on $G^{\mathbb C}/U$. 

The correspondence 
$\Lambda \mapsto F_\Lambda  : Z_{\frak t} \to H( G^{\mathbb C}/U, \underline{\mathbb C}^*)$ induces a
homomorphism. Also the correspondence $F \mapsto c_{1}(F)$ defines
a homomorphism of $H( G^{\mathbb C}/U, \underline{\mathbb C}^*)$ to $H^2(M, \mathbb Z)$. 
Then it is known that homomorphisms 
$Z_{\frak t} \overset{F}{\to} H( G^{\mathbb C}/U, \underline{\mathbb C}^*) \overset{c_1}{\to} H^2(M, \mathbb Z)
$
are in fact isomorphisms. In particular, the second Betti number $b_2(M)$ of $M$ is given by 
$b_2(M) = \dim {\frak t}$, the cardinality of  $\Pi - \Pi_0
$.

\subsection{K\"ahler-Einstein metrics on  a generalized flag manifold}
We set \newline $Z_{\frak t}^{+} = \left\{ \lambda \in Z_{\frak t} \ \big\vert \  (\lambda, \alpha) > 0 \ \  \mbox{for} \  \alpha \in \Pi - \Pi^{}_0 \right\}$. Then we have $\displaystyle Z_{\frak t}^{+} =  \sum_{\alpha \in  \Pi - \Pi^{}_0} {\mathbb Z}^+ \Lambda_{\alpha}$. 
 We define an element  $ \delta_{\frak m}  \in  \sqrt{-1}\frak h $  \, by \, 
 $\displaystyle  \delta_{\frak m} = \frac{1}{2} \sum_{\alpha \in \Delta_{\frak m}^{+} } \alpha$. 
 Let $c_1(M)$ be the first Chern class of $M$. Then we have that
 $ 2 \delta_{\frak m} \in Z_{\frak t}^{+}$ and  $c_1(M) = c_1(F_{2 \delta_{\frak m}})$.
 Put $ \displaystyle k_\alpha =  \frac{2 ( 2 \delta_{\frak m}, \alpha)}{(\alpha, \alpha)}$ for $ \alpha \in  \Pi - \Pi^{}_0$.  Then 
 $\displaystyle 2 \delta_{\frak m} =   \sum_{ \alpha  \in  \Pi - \Pi^{}_0}  k_\alpha \Lambda_\alpha
 = k_{\alpha_{i_1}}\Lambda_{\alpha_{i_1}} +\cdots + k_{\alpha_{i_r}}\Lambda_{\alpha_{i_r}} $ and each $k_{\alpha_{i_s}}$ is a positive integer. 
 The $G$-invariant metric $g_{ 2 \delta_{\frak m} }^{}$ on $G/K$ corresponding to $ 2 \delta_{\frak m} $, which is a K\"ahler Einstein metric, is given by 
$$ g_{ 2 \delta_{\frak m} }^{} =  \sum_{j_1, \dots, j_r} \left( \sum_{\ell = 1}^{r}
k_{\alpha_{i_{\ell}}}  j_{\ell} \frac{(\alpha_{j_\ell},  \alpha_{j_\ell})}{ 2} \right) B|_{ {\frak m}( j_1, \dots,  j_r )}.  
$$

\noindent
{\bf Example.} For the generalized flag manifold   $ G/K = Sp(n)/(U(p) \times U(q) \times Sp(n -p -q))$ where $n\ge 3, \ p, \, q \geq 1$,  we see that   $\Pi - \Pi^{}_0 = \{\alpha^{}_{p}, \  \alpha^{}_{p +q}\}$, where
$\alpha^{}_{p}$ and $\alpha^{}_{p +q}$ are two simple roots in the Dynkin diagram of $\frak{sp}(n)$ painted black.
 The painted Dynkinh diagram is shown below

\smallskip
\hspace{2.5cm}
\begin{picture}(160,45)(-40,-20)
\put(0, 0){\circle{4}}
\put(0,8.5){\makebox(0,0){$\alpha_1$}}
\put(0,-8){\makebox(0,0){2}}
\put(2, 0.3){\line(1,0){10}}
\put(20, 0){\makebox(0,0){$\ldots$}}
\put(28, 0.3){\line(1,0){10}}
\put(40, 0){\circle{4.4}}
\put(40,-8){\makebox(0,0){2}}
\put(40,8.5){\makebox(0,0){$\alpha_{p-1}$}}
\put(42, 0.3){\line(1,0){14}}
\put(58, 0){\circle*{4.4}}
\put(58,8.5){\makebox(0,0){$\alpha_p$}}
\put(58,-8){\makebox(0,0){2}}
\put(58, 0.3){\line(1,0){14}}
\put(74, 0){\circle{4.4}}
\put(74,-8){\makebox(0,0){2}}
\put(76, 0.3){\line(1,0){10}}
\put(94, 0){\makebox(0,0){$\ldots$}}
\put(102, 0.3){\line(1,0){10}}
\put(114, 0){\circle*{4}}
\put(114,8.5){\makebox(0,0){$\alpha_{p+q}$}}
\put(114,-8){\makebox(0,0){2}}
\put(116, 0.3){\line(1,0){10}}
\put(134, 0){\makebox(0,0){$\ldots$}}
\put(142, 0.3){\line(1,0){10}}
\put(154, 0){\circle{4}}
\put(154,8.5){\makebox(0,0){$\alpha_{n-1}$}}
\put(154,-8){\makebox(0,0){2}}
\put(159.2, 1.5){\line(1,0){13.8}}
\put(159.2, -1.2){\line(1,0){13.8}}
\put(155.46, -1.6){\scriptsize $<$}
\put(174.5, 0){\circle{4}}
\put(175,8.5){\makebox(0,0){$\alpha_n$}}
\put(174,-8){\makebox(0,0){1}}
\put(180,-1){\makebox(0,0){.}}  
\end{picture}

 It is $\displaystyle 2 \delta_{\frak m}  =  (p+q) \Lambda_{\alpha_{p}}  + (2 n -2 p - q +1 )\Lambda_{\alpha_{p+q}} $.  Then the K\"ahler-Einstein metric $g_{ 2 \delta_{\frak m} }^{}$ on $G/K$  is given by 
\begin{eqnarray} g_{ 2 \delta_{\frak m} }^{} =  (p+q) B|_{ {\frak m}( 1, 0 )} + (2 n - 2p -q +1) B|_{ {\frak m}( 0, 1 )}  \nonumber \\ + (2 n - p  +1) B|_{ {\frak m}( 1,  1)} + 2 (2 n - 2p -q +1) B|_{ {\frak m}( 0, 2 )}   \nonumber \\+(4 n-3p -q +2) B|_{ {\frak m}( 1, 2 )}+2(2 n - p + 1) B|_{ {\frak m}( 2, 2 )}.  \label{kahler_Einstein_metric}\end{eqnarray}

\subsection{Einstein metrics on  generalized flag manifolds with two isotropy summands}
We assume that the Lie group $G$ is simple. 
  For a generalized flag manifold $G/K$, we denote by $q$ the number of mutually non equivalent irreducible $\mbox{Ad}_G(K)$-modules  $ {\frak m}( j_1, \dots,  j_r )$ with 
 $\displaystyle \frak m = \sum_{j_1, \cdots, \, j_r} {\frak m}( j_1, \dots,  j_r ). 
  $
 
  If $q = 1$, it is known that $G/K$ is an irreducible Hermitian symmetric space with the symmetric pair $( \frak g, \frak k)$. 
 If $q = 2$, we have two $G$-invariant Einstein metrics on  $G/K$. 
    One is K\"ahler-Einstein metric and the other is non K\"ahler Einstein metric (\cite{ArCh1}).
 In fact, we see that the case  $q = 2$ occurs  only in the case $ r = b_2(G/K) =1$ and   $\displaystyle \frak m =  {\frak m}( 1 ) \oplus {\frak m}( 2 )$ (cf. \cite{ArChSa1}).
Note that   only $\displaystyle {2 \brack {11}} $ is non-zero. Put $ d_1= \dim {\frak m}( 1 ) $ and $ d_2 = \dim {\frak m}( 2 ) $.   For a $G$-invariant metric $  \langle  \,\,\, , \,\,\, \rangle = x_1\cdot B|_{{\frak m}( 1 ) } +  x_2\cdot B|_{{\frak m}( 2 )}$,   the components ${ r}_{1}, {r}_{2}$  of Ricci tensor ${r}$ of the metric $  \langle  \,\,\, , \,\,\, \rangle  $  are given by 
\begin{equation}
\left\{
\begin{array}{ll} 
r_1 &=  \displaystyle{\frac{1}{2 x_1} -
\frac{x_2}{2\,d_1\,{x_1}^2} {2 \brack 11}
\; }
\\ & \\
r_2 &=  \displaystyle{\frac{1}{2 x_2}  - \frac{1}{2\,d_2 \, x_2} {1 \brack 21}
 +
\frac{x_2}{4\,d_2\,{x_1}^2} {2 \brack 11}
\;.}
\end{array}
\right.
\end{equation}
 Since the metric $ ( \,\,, \,\,) = 1\cdot B|_{{\frak m}( 1 ) } +  2\cdot B|_{{\frak m}( 2 )}$ is K\"ahler Einstein, we see that $\displaystyle {2 \brack 11} = \frac{d_1  d_2}{d_1 + 4 d_2}$.
Note that a $G$-invariant metric $ \langle  \,\,\, , \,\,\, \rangle = x_1\cdot B|_{{\frak m}( 1 ) } +  x_2\cdot B|_{{\frak m}( 2 )}$ is Einstein if and only if $ r_1 = r_2$.  
 We normalize the equation $ r_1 = r_2$ by putting $ x_1 =1$. Then we see that the equation $ r_1 = r_2$ is reduced to a quadratic equation of $x_2$ and we have solutions $x_2 = 2$ and $\displaystyle x_2 = \frac{4 d_2}{d_1 + 2 d_2}$. 
 Since $\displaystyle x_2 = \frac{4 d_2}{d_1 + 2 d_2} \neq 2$,
  the Einstein metric  $\displaystyle 1\cdot B|_{{\frak m}( 1 ) } +  \frac{4 d_2}{d_1 + 2 d_2} \cdot B|_{{\frak m}( 2 )} $ is not K\"ahler.

 
 \subsection{Einstein metrics on  generalized flag manifolds with three isotropy summands}
 The case  $q = 3$ was studied by the author  and  Kimura  independently in \cite{Ar1} and \cite{Ki}.  
 We see that either  $ r = b_2(G/K) =1$ or  $ r = b_2(G/K) =2$. 
 
 For   the case of    $ r = b_2(G/K) =1$ we denote the $T$-roots system  of type $A_1(3)$,   that is  
  $\Delta_T^{+} = \{ \, \xi, 2 \xi,  3 \xi \, \}$.  There are seven cases. 
  The components ${ r}_{1}, {r}_{2}, {r}_{3}$  of the Ricci tensor ${r}$ of the metric 
$ \langle  \,\,\, , \,\,\, \rangle = x_1\cdot B|_{{\frak m}( 1 ) } +  x_2\cdot B|_{{\frak m}( 2 )}+
x_1\cdot B|_{{\frak m}( 3 ) }$  
  are given by
  \begin{equation*}
\left\{
\begin{array}{l} 
r_1  \displaystyle = \frac{1}{2 {x_1}}-\frac{ {x_2}}{2 {d_1}
   {x_1}^2} {2 \brack 11}+\frac{ 1}{2
   {d_1}}\left(\frac{{x_1}}{{x_2}
   {x_3}}-\frac{{x_3}}{{x_1}
   {x_2}}-\frac{{x_2}}{{x_1} {x_3}}\right) {3 \brack 12} \\
 r_2   \displaystyle = \frac{1}{2 {x_2}}+  \frac{1}{4 {d_2}} \left(\frac{x_2}{{x_1}^2}-\frac{2}{x_2}\right){2 \brack 11} 
  \displaystyle +\frac{1}{2 {d_2}} \left(\frac{x_2}{{x_1} {x_3}}-\frac{x_1}{{x_2}
   {x_3}}-\frac{x_3}{{x_1} {x_2}}\right) {3 \brack 12} \\
 r_3 \displaystyle  = \frac{1}{2 {x_3}}+  \frac{ 1}{2 {d_3}}
   \left( \frac{{x_3}}{{x_1}{x_2}} -\frac{{x_1}}{{x_2}{x_3}}-\frac{{x_2}}{{x_1} {x_3}}\right) {3 \brack 12}. 
 \end{array}
\right.
\end{equation*}
The system of equations $r_1 = r_2 = r_3 $ reduces to a polynomial equation of degree 5.
There is a unique K\"ahler-Einstein metric and two non K\"ahler Einstein metrics. These were
explicitely found  in \cite{AnCh}.

For the  case $q = 3$ and $b_2(G/K)  =2$ 
we denote  the case of    $ r = b_2(G/K) =2$ that $T$-roots system is of type $A_2$,  that is
  $\Delta_T^{+} = \{ \, \xi_1,   \xi_2,    \xi_1 + \xi_2 \, \}$.  
There are three cases here. 
The space $\SU(n)/\s(\U(\ell)\!\times\! \U(m)\!\times\! \U(k))$ ($ n=\ell+m+k$)
admits three K\"ahler-Einstein metrics (for  $\ell, m, l$ distinct) and one non K\"ahler Einstein metric.
The space $\SO(2n)/(\U(n-1)\!\times\! \U(1)) $ admits two K\"ahler-Einstein metrics and one non-K\"ahler Einstein metric.
Finally, the space $E_6/(\SO(8) \!\times\! \U(1)\!\times\! \U(1))$ admits
a unique K\"ahler-Einstein and a unique non K\"ahler Einstein metric.  This metric 
is normal.

\subsection{Einstein metrics on  generalized flag manifolds with four isotropy summands}
We see that in this case either  $ r = b_2(G/K) =1$ or  $ r = b_2(G/K) =2$.
For the case $ r = b_2(G/K) =1$ we denote the $T$-root system of type $A_1(4)$,  that is  
  $\Delta_T^{+} = \{ \, \xi, 2 \xi,  3 \xi, \, 4 \xi \, \}$. 
  There are four  flag manifolds of this type and   $G$ is always  an exceptional Lie group.
  It follows that they admit one K\"ahler-Einstein and two non K\"ahler Einstein metrics.
  
   For the case of   $ r = b_2(G/K) = 2$ we denote the $T$-root system  of type $B_2$,  that is  
  $\Delta_T^{+} = \{ \, \xi_1,   \xi_2,    \xi_1 + \xi_2, \,   2\xi_1 + \xi_2 \, \}$. There are six
  flag manifolds of this type.
   We divide the cases of type $B_2$
   into   $B_2(a)$ and $B_2(b)$. 
 
 The case of $B_2(a)$ is  one of $\SO(2 \ell +1)/ \SO(2 \ell -3) \times \U(1) \times \U(1)$, $\SO(2 \ell)/ \SO(2 \ell -4) \times \U(1) \times \U(1)$,  $E_6/ \SU(5)\times \U(1) \times \U(1)$ or $E_7/ \SO(10)\times \U(1) \times \U(1)$. 
 The case of $B_2(b)$ is  one of $\Sp( \ell )/ \U(p) \times \U(\ell-p)$ or  $\SO(2 \ell)/  \U(p) \times \U(\ell-p)$. These results are presented in \cite{ArCh2}.
 
 For  the case of $B_2(a)$, there exist exactly eight $G$-invariant Einstein metrics on $G/K$, four  of them are K\"ahler-Einstein and the others are non-K\"ahler Einstein.
All above results are presented in \cite{ArCh2}.
 
 The case of  $B_2(b)$ is the most difficult one.  
 For $\SO(2 \ell)/ ( \U(p) \times \U(\ell-p) )$ ( $\ell \geq 4 $ and $2 \leq p \leq \ell -2$), there exist four non-K\"ahler Einstein metrics  for the pairs $(\ell, p)$ $= (12, 6)$, $ (10, 5)$, $ (8, 4)$, $ (7, 4)$,  $ (7, 3)$,  $ (6, 4)$, $ (6, 3)$,  $ (6, 2)$,  $ (5, 3)$, $ (5, 3)$,  $ (4, 2)$.  For the other pairs  $(\ell, p)$, there exit two non-K\"ahler Einstein metrics.
 The complete study if this space
 as well as the isometry problem is presented in the paper \cite{ArChSa1}.
 
 Finally, the space $\Sp( \ell )/ (\U(p) \times \U(\ell-p) )$ admits four isometric
 K\"ahler-Einstein metrics and two non K\"ahler Einstein metrics (cf. \cite{ArChSa1}).  These metrics are isometric, as it was
 proved in \cite{ArChSa5}.

 In order to get an idea of the difficulty of such parametric systems of equations, we mention that  the Einstein equation for $\Sp(n)/(\U(p)\times \U(n-p))$ reduces to the system
  $r_1-r_3=0, r_1-r_2=0,  r_3-r_4=0$, which is equivalent to the parametric system  
  
  \begin{equation*}
 \left.
 \begin{tabular}{l}
$(x_1-x_3)(x_1 x_2 + p x_1 x_2 + x_2 x_3 + p x_2 x_3 + x_1 x_4   + n x_1 x_4$\\
 $-  p x_1 x_4 - 2 x_2 x_4    - 2 n x_2 x_4 + x_3 x_4 + n x_3 x_4 - p x_3 x_4)=0,$\\
$4(n+1) x_3 x_4 (x_2-x_1)+ (n+p+1)x_4(x_1^2-x_2^2)- (n-3p+1)x_3^2x_4$\\
$ +(p+1)x_2(x_1^2-x_3^2-x_4^2)=0,$\\
$4(n+1)x_1x_2(x_4-x_3)+ (2n-p+1)x_2(x_3^2-x_4^2)+ (2n-3p-1)x_1^2x_2$\\
$ +(n-p+1)x_4(x_3^2-x_1^2-x_2^2)=0$\\
\end{tabular}
\right\} 
  \end{equation*} 
  
  \noindent
  for the unknowns $x_1, x_2, x_3, x_4>0$.
  
 \subsection{Einstein metrics on  generalized flag manifolds with five isotropy summands}
 The classification of generalized flag manifolds $M=G/K$ with five isotropy summands  was obtained in \cite{ArChSa4}.
They are obtained in the following ways: paint black one  simple root of Dynkin mark 5, that is 
 $\Pi\backslash\Pi_{0}=\{\al_{p} : \Mark(\al_{p})=5\}$, or 
  paint black two simple roots, one of Dynkin mark 1 and  one of Dynkin mark 2,   that is $\Pi\backslash\Pi_{0}=\{\al_{i}, \al_{j} : \Mark(\al_{i})=1, \ \Mark(\al_{j})=2\}$, or  
  paint black two simple roots both of Dynkin mark 2,  that is $\Pi\backslash\Pi_{0}=\{\al_{i}, \al_{j} : \Mark(\al_{i})= \Mark(\al_{j})=2\}$.  
 They correspond to the spaces $E_8/\U(1)\times \SU(4)\times\SU(5)$, $\SO(2\ell+1)/\U(1)\times\U(p)\times
 \SO(2(\ell-p-1)+1)$, $\SO(2\ell)/\U(1)\times\U(p)\times
 \SO(2(\ell-p-1))$, $E_6/\SU(4)\times \SU(2)\U(1)\times \U(1)$ and $E_6/\SU(6)\times \U(1)\times \U(1)$.  
 
The main difficulty in constructing the  Einstein equation for a $G$-invariant metric on one of the above homogeneous spaces  is  the calculation of the non zero structure constants $\displaystyle{k \brack {ij}}$ of $G/K$ with respect to the decomposition 
$T_{o}M\cong\fr{m}=\fr{m}_1\oplus\fr{m}_2\oplus\fr{m}_3\oplus\fr{m}_4\oplus\fr{m}_5$ of the tangent space
of $M$. 
A first step towards this procedure, is to use the known K\"ahler-Einstein metric (cf. Section 4).
Secondly, and this is was the main contribution of the papers \cite{ArChSa4} and \cite{ChSa},   we can take advantage of the fibration of a flag manifold  over another  flag manifold,  and  by using methods of   Riemannian submersions it possible to compare the Ricci tensors of the total and  base spaces, respectively.
Such fibration is a an extension of the well known twistor fibration of a generalized flag manifold over a symmetric space (cf. \cite{BuRa}).
In this way we can calculate  $\displaystyle{k \brack {ij}}$ in terms of the dimension of the associated submodules $\fr{m}_i$.

 The Einstein equation  reduces to a  polynomial system of four equations in four unknowns.
  For the exceptional flag manifolds it is possible to  classify all homogeneous Einstein metrics.    
 For the classical flag manifolds a complete classification of homogeneous Einstein metrics in the general case  is a difficult task, because the corresponding systems of equations depend on  four positive parameters (which define the invariant Riemannian metric), the Einstein constant $\lambda>0$ and the positive integers $\ell$ and $p$.    
 However, by using Gr\"obner bases we can show that the equations are reduced to a polynomial equation of   one variable and then we can prove  the existence of non K\"aler Einstein metrics.  
 In fact, this is another contribution of \cite{ArChSa4}, because we  show existence of real solutions
 for polynomial equations whose coefficients depend on parameters ($\ell$ and $p$).  Finally the isometry question for the Einstein metrics found, is also answered in this work.
 The results of \cite{ArChSa4} and \cite{ChSa} can be summarized as follows:
 
\begin{theorem}
 Let $M=G/K$ be one of the flag manifolds $\E_6/(\SU(4)\times \SU(2)\times\U (1)\times\U (1))$ or
 $\E_7/(\U(1)\times\U(6))$. Then $M$ admits exactly seven $G$-invariant Einstein metrics up to isometry.
 There are two K\"ahler-Einstein metrics and five non K\"ahler Einstein metrics (up to scalar).
 \end{theorem}
 
\begin{theorem} Let $M=G/K$ be one of the flag manifolds $\SO(2\ell +1)/(\U(1)\times\U (p)\times\SO(2(\ell -p-1)+1))$  $(\ell \ge 3, \ 3\le p\le\ell -1)$ or    
 $\SO(2\ell)/(\U(1)\times\U (p)\times\SO(2(\ell -p-1)))$ $(\ell\ge 5, \ 3\le p\le\ell -3)$.
 Then $M$ admits at least two $G$-invariant non K\"ahler Einstein metrics.  
 \end{theorem}
 
\begin{theorem} Let $M=G/K$ be one of the flag manifolds 
 $\SO(2\ell +1)/(\U(1)\times\U (2)\times\SO(2\ell -5))$  $(\ell \ge 6)$ or    
 $\SO(2\ell)/(\U(1)\times\U (2)\times\SO(2(\ell -3)))$ $(\ell\ge 7)$.
 Then $M$ admits at least four $G$-invariant non K\"ahler Einstein metrics.  
 \end{theorem}
 
\begin{theorem} The flag manifold $\E_8/\U(1)\times \SU(4)\times\SU (5)\times\U (1)$ 
 admits exactly six $E_8$-invariant Einstein metrics up to isometry.
 One is K\"ahler-Einstein metric and five are non K\"ahler Einstein metrics (up to scalar).
 \end{theorem}

  \subsection{Einstein metrics on  generalized flag manifolds with six isotropy summands}
Invariant Einstein metrics on generalized flag manifolds with six isotropy summands have not been completely classified.  
  A typical example is the flag manifold $G_2/T$, where $T$ is a maximal torus in $G_2$. The painted Dynking diamgram is shown below
  
  \hspace{5cm}
  \begin{picture}
 (85,45)(25,-25)
\put(52, 0){\circle*{4}}
\put(52,8.5){\makebox(0,0){$\alpha_{1}$}}
\put(52,-8){\makebox(0,0){3}}
\put(53., 1.8){\line(1,0){15}}
\put(53., -1.8){\line(1,0){15}}
\put(54., 0){\line(1,0){18}}
\put(62, -2.448){ $>$}
\put(74.3, 0){\circle*{4}}
\put(74,8.5){\makebox(0,0){$\alpha_2$}}
\put(73,-8){\makebox(0,0){2}}
\end{picture}

     The highest root   $\widetilde\alpha$  of ${\frak g_2 ^{\mathbb C}}$ is given by $\widetilde\alpha = 2 \alpha_1+ 3\alpha_2$. Thus we have a decomposition of $\frak m$ into  six mutually non-equivalent irreducible $\mbox{Ad}_G(H)$-modules  $ {\frak m}(j_1,  j _2)$ :
  $\displaystyle \frak m =  {\frak m}( 1, 0 ) \oplus {\frak m}( 0, 1 ) \oplus {\frak m}( 1, 1 ) \oplus {\frak m}( 1, 2 ) \oplus {\frak m}( 1, 3 ) \oplus {\frak m}( 2,3 ). 
  $ 
   There is only one  complex strucure and thus,  up to isometry, there exist only one  K\"ahler-Einstein metric. 
  There exist two non-K\"ahler Einstein metrics up to isometry.  These are obtained from solutions of polynomial of degree 14 (cf. \cite{ArChSa6}). 
There are four other generalized flag manifolds (all determined by the exceptional Lie groups, $F_4$, $E_6$, $E_7$, $E_8$) with   $T$-roots being of $G_2$ type.   For all cases, there is  only one  K\"ahler-Einstein metric and  at least six non K\"ahler Einstein metrics up to isometry. 

Other examples are the   flag manifold $\SU(4)/T$ and $\SU(10)/\s(\U(1)\times \U(2)\times \U(3)\times \U(4))$.
Note that for these cases  $q =6$ and the system of $T$-roots is of type $A_3$. 
 For the case $\SU(4)/T$   there is only one  complex strucure and thus,  up to isometry, there exists only one  K\"ahler-Einstein metric. There exist three non K\"ahler Einstein metrics up to isometry, one of them is normal (cf. \cite{Sa}).
  For the case $\SU(10)/\s(\U(1)\times \U(2)\times \U(3)\times \U(4))$ there are twelve  complex strucure and thus,  up to isometry, there exist twelve  K\"ahler-Einstein metrics. 
  There exist twelve non K\"ahler Einstein metrics up to isometry.  These are obtained from solutions of polynomial of degree 68.   
  
  In \cite{ArChSa2} all invariant Einstein metrics were found for the generalized flag manifolds $\Sp(3)/(\U(1) \times \U(1) \times \Sp(1))$, $\Sp(4)/(\U(1) \times \U(1) \times \Sp(2))$ and $\Sp(4)/(\U(2) \times \U(1) \times \Sp(1))$, and in \cite{ArChSa3} generalized flag manifolds  with $G_2$-type $\frak t$-roots (apart from  the full flag manifold $G_2/T$), were classified.  The main result is the following:
  
  \begin{theorem}\label{theorem1}
A generalized flag manifold $G/K$ with $G_2$-type $\frak t$-roots, which is not the full flag manifold $G_2/T$, admits exactly  one  invariant K\"ahler Einstein metric  and  six non K\"ahler invariant Einstein metrics up to isometry and scalar. 
These are the spaces $F_4/ \U(3)\times \U(1)$, $E_6/ \U(3)\times  \U(3)$, $ E_7/ \U(6)\times  \U(1)$ and $E_8/ E_6\times  \U(1)\times  \U(1)$. 
The  full flag manifold $G_2/T$ admits exactly  one  invariant K\"ahler Einstein metric  and  two non K\"ahler invariant Einstein metrics up to isometry and scalar.
\end{theorem}
Also, in \cite{ChSa} all $E_8$-invariant Einstein metrics on $E_8/\U(1)\times\SU(2)\times\SU(3)\times\SU(5)$ were classified and in \cite{WaZh} all invariant Einstein metrics were classified on
$F_4/\U(1)\times\U(1)\SO(5)$, $E_6/\U(1)\times\U(1)\times\SU(2)\times\SU(3)\times\SU(2)$ and
$E_8/\SO(12)\times\U(1)\times\U(1)$.

\section{Homogeneous spaces with equivalent isotropy summands}

If the isotropy representation of a homogeneous space $G/K$ 
$
{\frak m} = {\frak m}_1 \oplus \cdots \oplus {\frak m}_q$ contains equivalent summands, then the description
of all $G$-invariant metrics is more complicated, hence the problem of finding all invariant Einstein metrics is quite difficult.  We refer to \cite{St} for a study of invariant metrics on homogeneous spaces with equivalent isotropy summands.

An important class of homogeneous spaces with such a property are the real Steifel manifolds $V_k\mathbb{R} ^n=\SO(n)/\SO(n-k)$ of orthonormal $k$-frames in $\bb{R} ^n$ (as well as the complex and quaternionic Stiefel manifolds $V_k\mathbb{C} ^n=\SU(n)/\SU(n-k)$ and $V_k\mathbb{H} ^n=\Sp(n)/\Sp(n-k)$ of orthonormal $k$-frames in $\bb{C} ^n$ and $\bb{H} ^n$ respectively). 
 The simplest case is
   the sphere $\mathbb{S}^{n-1}=\SO(n)/\SO(n-1)$, which is an irreducible symmetric space, therefore it admits up to scale a unique invariant Einstein metric.
Concerning Einstein metrics on other Stiefel manifolds we recall the following:
   In \cite{K} S. Kobayashi proved existence of an invariant Einstein metric on the unit tangent bundle $T_1S^n = \SO(n)/\SO(n-2)$.  In \cite{S} A. Sagle proved that the Stiefel manifolds 
   $V_k\mathbb{R} ^n=\SO(n)/\SO(n-k)$ admit at least one homogeneous Einstein metric.  
   For $k\ge 3$ G. Jensen in \cite{J2} found a second metric.
  For $n=3$ the Lie group $\SO(3)$ admits a unique Einstein metric.
  For $n\ge 5$  A. Back and W.Y. Hsiang in \cite{BH} proved that $\SO(n)/\SO(n-2)$ admits
  exactly one homogeneous Einstein metric. The same result was obtained by M. Kerr in \cite{Ke} by proving that the diagonal metrics are the only invariant metrics on $V_2\mathbb{R}^n$ (see also \cite{A1}, \cite{AK}).
  The Stiefel manifold $\SO(4)/\SO(2)$ admits exactly two invariant Einstein metrics (\cite{ADF}).
  One is Jensen's metric and the other one is the product metric on $S^3\times S^2$.
  
  Finally, in \cite{ADN1} the  author, V.V. Dzhepko and Yu. G. Nikonorov proved that
  for $s>1$ and $\ell >k\ge 3$ the Stiefel manifold $\SO(s k+\ell)/\SO(\ell)$ admits at least four
  $\SO(s k+\ell)
  $-invariant Einstein metrics, two of which are Jensen's metrics.
  The special case $\SO( 2 k + \ell)/\SO(\ell)$ admitting at least four
  $\SO(2 k+\ell) $-invariant Einstein metrics was treated in \cite{ADN2}.
  Corresponding results for the quaternionic Stiefel manifolds $\Sp(sk+\ell)/\Sp(\ell)$ were obtained in
  \cite{ADN3}.
  
  In the recent work \cite{ASS} it was shown that
  the Stiefel manifold $V_4\mathbb{R}^n=\SO(n)/\SO(n-4)$ admits two more $\SO(n)$-invariant Einstein metrics and that
  $\SO(7)/\SO(2)$ admits four more $\SO(7)$-invariant Einstein metrics  (in addition to the ones obtained in \cite{ADN1}).  This was achieved  by making appropriate symmetry assumptions on the set of $\SO(n)$-invariant metrics on $V_4\mathbb{R}^n$.  We refer to \cite{ADN1} and \cite{St} for a  more general presentation of such method.

\subsection{The Stiefel manifolds $V_{k_1+k_2}\mathbb{R}^{k_1+k_2+k_3} = \SO(k_1+k_2+k_3)/\SO(k_3)$}

We first consider the homogeneous space $G/K= \SO(k_1+k_2+k_3)/(\SO(k_1)\times\SO(k_2)\times\SO(k_3))$, where the embedding of $K$ in $G$ is diagonal. This is an example of a generalized Wallach space.  These spaces have been recently classified by Yu.G. Nikonorov (\cite{Ni2}) and Z. Cheng, Y. Kang and K. Liang (\cite{ChKaLi}).
If $\lambda _n$ denotes the standard representation of $\SO(n)$, then $\Ad^{\SO(n)}=\wedge ^2\lambda _n$.
Let $\sigma _i:\SO(k_1)\times\SO(k_2)\times\SO(k_3)\to \SO(k_i)$  be the projection onto
the factor $\SO(k_i)$ ($i=1,2,3$) and let $p_{k_i}=\lambda _{k_i}\circ\sigma _i$.
Then a calculation shows that 
$$
\left.\Ad ^G\right |_K=\Ad ^K\oplus (p_{k_1}\otimes p_{k_2})\oplus (p_{k_1}\otimes p_{k_3})
\oplus (p_{k_2}\otimes p_{k_3}).
$$
Thus the isotropy representation of $G/K$ is
$(p_{k_1}\otimes p_{k_2})\oplus (p_{k_1}\otimes p_{k_3})
\oplus (p_{k_2}\otimes p_{k_3})$
and this induces a decomposition of 
the tangent space $\fr{m} $ of $G/K$  into  three  $\Ad(K)$-submodules 
$$
\fr{m} = \fr{m}_{12}\oplus  \fr{m}_{13}\oplus  \fr{m}_{23}. 
$$
In fact, $\fr{m}$ is given by  $\fr{k}^{\perp} $ in $ \fr{g} = \fr{so}(k_1+ k_2+k_3)$ with respect to  $-B$. If we denote by $M(p,q)$ the set of all $p \times q$ matrices, then we see that $\fr{m}$ is given by 
\begin{equation*}
\fr{m}=  \left\{\begin{pmatrix}
 0 & {A}_{12} & {A}_{13}\\
 -{}^{t}_{}\!{A}_{12} & 0 & {A}_{23}\\
 -{}^{t}_{}\!{A}_{13} & -{}^{t}_{}\!{A}_{23} & 0 
 \end{pmatrix} \  \Big\vert \  {A}_{12} \in M(k_1, k_2), {A}_{13} \in M(k_1, k_3), {A}_{23} \in M(k_2, k_3) \right\} 
\end{equation*}
and we have 
\begin{equation*}
 \fr{m}_{12}= \begin{pmatrix}
 0 & {A}_{12} & 0\\
 -{}^{t}_{}\!{A}_{12} & 0 &0\\
0  & 0 & 0 
 \end{pmatrix},  \quad  
 \fr{m}_{13}= \begin{pmatrix}
 0 & 0 &{A}_{13}\\
0 & 0 &0\\
 -{}^{t}_{}\!{A}_{13}   & 0 & 0 
 \end{pmatrix},  \quad  
 \fr{m}_{23}= \begin{pmatrix}
 0 & 0 & 0\\
0 & 0 &{A}_{23}\\
0  &  -{}^{t}_{}\!{A}_{23} & 0 
 \end{pmatrix}.  
\end{equation*}
Note that the action of $\Ad(k)$ ($k \in K$) on $ \fr{m}$ is given by 
\begin{equation*}
\Ad(k) \begin{pmatrix}
 0 & {A}_{12} & {A}_{13}\\
 -{}^{t}_{}\!{A}_{12} & 0 & {A}_{23}\\
 -{}^{t}_{}\!{A}_{13} & -{}^{t}_{}\!{A}_{23} & 0 
 \end{pmatrix}  = 
  \begin{pmatrix}
 0 &{}^t h_1 {A}_{12} h_2 & {}^t h_1{A}_{13} h_3\\
 -{}^{t}_{}h_2 {}^{t}_{}\!{A}_{12}h_1 & 0 & {}^t h_2{A}_{23} h_3\\
 -{}^{t}_{}h_3{}^{t}_{}\!{A}_{13} h_1 & -{}^{t}_{}h_3{}^{t}_{}\!{A}_{23} h_2& 0 
 \end{pmatrix}, 
\end{equation*}
 where $ \begin{pmatrix}
 h_1 & 0& 0\\
0 & h_2 &0\\
0  & 0 & h_3 
 \end{pmatrix} \in K$.  Thus  the irreducible submodules  $\fr{m}_{12}$,  $\fr{m}_{13}$ and  $\fr{m}_{23}$  are  mutually non equivalent. 
 
We now consider the Stiefel manifold $G/H=\SO(k_1+k_2+k_3)/\SO(k_3)$ and we take into account the diffeomorphism
$
G/H=(G\times \SO(k_1)\times\SO(k_2))/((\SO(k_1)\times\SO(k_2))\times\SO(k_3))=\widetilde{G}/\widetilde{H}.
$
Let
$
\fr{p}=\fr{so}(k_1)\oplus\fr{so}(k_2)\oplus  \fr{m}_{12}\oplus  \fr{m}_{13}\oplus  \fr{m}_{23} 
$
be an $\Ad(\SO(k_1)\times\SO(k_2)\times\SO(k_3))$-invariant decomposition of the tangent space $\fr{p}$ of $G/H$ at $eH$, where the corresponding submodules are non equivalent.  Then we consider a subset of all $G$-invariant metrics on $G/H$  determined by the $\Ad(\SO(k_1)\times\SO(k_2)\times\SO(k_3))$-invariant scalar products on $\fr{p}$ given by
\begin{equation}\label{metric1}
\langle \   ,\  \rangle =  x_1 \, (-B) |_{\fr{so}(k_1)}+ x_2 \, (-B) |_{ \fr{so}(k_2)}
 + x_{12} \,  (-B) |_{ \fr{m}_{12}}+ x_{13} \,  (-B) |_{ \fr{m}_{13}} + x_{23} \,  (-B) |_{ \fr{m}_{23}},
\end{equation}
 for $k_1 \geq 2$, $k_2 \geq 2$ and $k_3 \geq 1$.  Since the submodules are non equivalent we can use Lemma
 \cite{PS} and obtain the following:
 
\begin{lemma}\label{lemma5.4}
The components  of  the Ricci tensor ${r}$ for the invariant metric $\langle \   ,\  \rangle$ on $G/H$ defined by  {\rm (\ref{metric1})} are given as follows:  
\begin{equation}\label{eq18}
\left. 
\small{\begin{array}{lll} 
r_1 &=&  \displaystyle{\frac{k_1-2}{4 (n -2)  x_1} +
\frac{1}{4 (n -2) } \biggl(k_2 \frac{x_1}{{x_{12}}^2}} +k_3 \frac{x_1}{{x_{13}}^2} \biggr), 
 \\   \\
 r_2 &=&  
\displaystyle{\frac{k_2-2}{4 (n -2)  x_2} +
\frac{1}{4 (n -2)} \biggl(k_1 \frac{x_2}{{x_{12}}^2} +k_3 \frac{x_2}{{x_{23}}^2} \biggr),} 
\\  \\
r_{12} &=&   \displaystyle{\frac{1}{ 2 x_{12}} +\frac{k_3}{4 (n -2)}\biggl(\frac{x_{12}}{x_{13} x_{23}} - \frac{x_{13}}{x_{12} x_{23}} - \frac{x_{23}}{x_{12} x_{13}}\biggr) }\\  \\ & &
\displaystyle{-
\frac{1}{4 (n -2)} \biggl( (k_1-1) \frac{x_1}{{x_{12}}^2} + (k_2-1) \frac{x_2}{{x_{12}}^2} \biggr)},
\\  \\
r_{13}  &=&   \displaystyle{\frac{1}{  2 x_{13}} +\frac{k_2}{4 (n -2)}\biggl(\frac{x_{13}}{x_{12} x_{23}} - \frac{x_{12}}{x_{13} x_{23}} - \frac{x_{23}}{x_{12} x_{13}}\biggr) -
\frac{1}{4 (n -2)} \biggl( (k_1-1) \frac{x_1}{{x_{13}}^2}  \biggr)}
\\ \\
r_{23}  &=&   \displaystyle{\frac{1}{ 2 x_{23}} +\frac{k_1}{4 (n -2)}\biggl(\frac{x_{23}}{x_{13} x_{12}} - \frac{x_{13}}{x_{12} x_{23}} - \frac{x_{12}}{x_{23} x_{13}}\biggr) -
\frac{1}{4 (n -2)} \biggl(  (k_2-1) \frac{x_2}{{x_{23}}^2} \biggr)}.   
\end{array} } \right\}
\end{equation}
where $n = k_1+k_2+k_3$. 
\end{lemma}

For  $k_1 =1$  we have  the Stiefel manifold $G/H=\SO(1+k_2+k_3)/\SO(k_3)$ with corresponding decomposition 
$
\fr{p} =  \fr{so}(k_2)\oplus  \fr{m}_{12} \oplus  \fr{m}_{13} \oplus  \fr{m}_{23}. 
$
We then consider $G$-invariant metrics on $G/H$ determined by the  $\Ad(\SO(k_2)\times\SO(k_3))$-invariant scalar products on $\fr{p}$ given by 
\begin{equation} \label{metric2} 
\langle \   ,\ \rangle =     x_2 \, (-B) |_{ \fr{so}(k_2)} + x_{12} \,  (-B) |_{ \fr{m}_{12}} + x_{13} \,  (-B) |_{ \fr{m}_{13}} + x_{23} \,  (-B) |_{ \fr{m}_{23}}
\end{equation}

\begin{lemma}\label{metric22}
The components  of  the Ricci tensor ${r}$ for the invariant metric $ \langle \  \ ,\ \ \rangle $ on $G/H$ defined by  \em{(\ref{metric2})}, are given as follows:  
\begin{equation}\label{eq19}
\left. {\small \begin{array}{l} 
 r_2  = 
\displaystyle{\frac{k_2-2}{4 (n -2) x_2} +
\frac{1}{4 (n -2)} \biggl(  \frac{x_2}{{x_{12}}^2} +k_3 \frac{x_2}{{x_{23}}^2} \biggr)},
\\  \\
r_{12}   =  \displaystyle{\frac{1}{2 x_{12}} +\frac{k_3}{4 (n -2)}\biggl(\frac{x_{12}}{x_{13} x_{23}} - \frac{x_{13}}{x_{12} x_{23}} - \frac{x_{23}}{x_{12} x_{13}}\biggr) -
\frac{1}{4 (n -2)} \biggl(  (k_2-1) \frac{x_2}{{x_{12}}^2} \biggr)},
\\  \\
r_{23}   =  \displaystyle{\frac{1}{2 x_{23}} +\frac{1}{4 (n -2)}\biggl(\frac{x_{23}}{x_{13} x_{12}} - \frac{x_{13}}{x_{12} x_{23}} - \frac{x_{12}}{x_{23} x_{13}}\biggr) -
\frac{1}{4 (n -2)} \biggl(  (k_2-1) \frac{x_2}{{x_{23}}^2} \biggr)}, 
\\  \\
r_{13}   =  \displaystyle{\frac{1}{ 2 x_{13}} +\frac{k_2}{4 (n -2)}\biggl(\frac{x_{13}}{x_{12} x_{23}} - \frac{x_{12}}{x_{13} x_{23}} - \frac{x_{23}}{x_{12} x_{13}}\biggr) }. 
\end{array} } \right\} 
\end{equation}
where $n = 1+k_2+k_3$. 
\end{lemma}

\subsection{The Stiefel manifold $V_5\mathbb{R} ^n$}

For the Stiefel manifold $V_5\mathbb{R} ^n = \SO(n)/\SO(n-5)$ we can let $k_1=2, k_2=3, k_3=n-5$ and consider  $\Ad(\SO(2)\times\SO(3)\times\SO(n-5))$-invariant scalar products of the form (\ref{metric1}), or let $k_1=1, k_2=4, k_3=n-5$ and consider $\Ad(\SO(4)\times\SO(n-5))$-invariant scalar products of the form (\ref{metric2}).

\begin{theorem}\label{theorem_v4_R^n} 
The Stiefel manifold $V_5\mathbb{R}^n = \SO(n)/\SO(n-5)$ ($n\ge 7$)  admits at least six invariant Einstein metrics. Two of them are Jensen's  metrics, two are given by  $\Ad(\SO(4)\times\SO(n-5))$-invariant scalar products of the form (\ref{metric2}), and the other two are given by  $\Ad(\SO(2)\times\SO(3)\times\SO(n-5))$-invariant scalar products of the form (\ref{metric1}).  
\end{theorem} 
\begin{proof}


We consider $\Ad(\SO(4)\times\SO(n-5))$-invariant inner products of the form \ref{metric2}. 
Then from Lemma \ref{metric22} we see that the components of the Ricci tensor $r$ are given by:
\begin{equation}\label{eq10}
\left. {\small \begin{array}{l} 
 r_2  = \displaystyle{\frac{1}{4(n-2) x_2} + \frac{1}{4(n-2)} \biggl(\frac{x_2}{{x_{12}}^2} +(n-5) \frac{x_2}{x_{23}^2}\biggr)},
\\  \\
r_{12}   =  \displaystyle{\frac{1}{2 x_{12}} +\frac{n-4}{4(n-2)}\biggl(\frac{x_{12}}{x_{13} x_{23}} - \frac{x_{13}}{x_{12} x_{23}} - \frac{x_{23}}{x_{12} x_{13}}\biggr) -
\frac{3}{4(n-2)}  \frac{x_2}{{x_{12}}^2}},
\\  \\
r_{23}   =  \displaystyle{\frac{1}{2 x_{23}} +\frac{1}{4(n-2)}\biggl(\frac{x_{23}}{x_{13} x_{12}} - \frac{x_{13}}{x_{12} x_{23}} - \frac{x_{12}}{x_{23} x_{13}}\biggr) -
\frac{3}{4(n-2)} \frac{x_2}{{x_{23}}^2}},
\\  \\
r_{13}   =  \displaystyle{\frac{1}{ 2 x_{13}} +\frac{1}{(n-2)}\biggl(\frac{x_{13}}{x_{12} x_{23}} - \frac{x_{12}}{x_{13} x_{23}} - \frac{x_{23}}{x_{12} x_{13}}\biggr) }.
\end{array} } \right\} 
\end{equation}

We consider the system of equations   
 
\begin{equation}\label{eq21a} 
r_2 =  r_{12}, \, \,  r_{12} = r_{23}, \, \,  r_{23} = r_{13}. 
\end{equation}
We put $x_{23} = 1$ and from  system (\ref{eq21a}) we have:
\begin{equation}\label{eq11}
\left.{ 
\begin{array}{lll}
f_{1} &=& -n x_{2} x_{12}^3+5 x_{2} x_{12}^3+n x_{13} x_{2}^2 x_{12}^2-5
   x_{13} x_{2}^2 x_{12}^2+2 x_{13} x_{12}^2+n x_{13}^2 x_{2}x_{12}\\
   && -5 x_{13}^2 x_{2} x_{12}+n x_{2} x_{12}-2 n x_{13} x_{2}
   x_{12}+4 x_{13} x_{2} x_{12}-5 x_{2} x_{12}+4 x_{13} x_{2}^2=0\\
f_{2} &=& n x_{12}^3-4 x_{12}^3-2 n x_{13} x_{12}^2+4 x_{13} x_{12}^2+3 x_{13}
   x_{2} x_{12}^2-n x_{13}^2 x_{12}+6 x_{13}^2 x_{12}\\
  && -n x_{12}+2 n x_{13} x_{12}-4 x_{13} x_{12}+4 x_{12}-3 x_{13} x_{2}=0  \\
f_{3} &=& 3 x_{12}^2-2 n x_{12}+2 n x_{13} x_{12}-4 x_{13} x_{12}-3 x_{13} x_{2} x_{12}+4 x_{12}-5 x_{13}^2+5 =0.   
\end{array} } \right\} 
\end{equation}
We consider a polynomial ring $R= {\mathbb Q}[z,  x_2, x_{12}, x_{13}] $ and an ideal $I$ generated by 
$\{ f_1, \, f_2, \, f_3, $  $  \,z  \, x_2 \, x_{12} \, x_{13} -1\}  
$  to find non zero solutions of the above equations. 
We take a lexicographic order $>$  with $ z >   x_2 > x_{12} > x_{13}$ for a monomial ordering on $R$. Then by the aid of computer, we see that a  Gr\"obner basis for the ideal $I$ contains the  polynomial
$(x_{13} - 1) \, h_{1}(x_{13}),$
where $h_{1}(x_{13})$ is a polynomial of $x_{13}$ given by 
{\small \small\small  \begin{equation*} 
\begin{array}{l}
h_{1}(x_{13})= (-1 + n)^3 (-7 + 3 n)^2 (-81 + 78 n - 22 n^2 + 2 n^3) (-1 + 22 n - 
    14 n^2 + 2 n^3)x_{13}^{10} \\
 -2 (-1 + n)^2 (-7 + 3 n) (73599 - 286463 n + 441002 n^2 - 361584 n^3 + 175526 n^4 - 52216 n^5 \\ 
 + 9376 n^6 - 936 n^7 + 40 n^8)x_{13}^{9} + (-1 + n) (11395069 - 47604434 n + 87540578 n^2\\ 
 - 93088081 n^3 +63294239 n^4 - 28742935 n^5 + 8835132 n^6 - 1817860 n^7 \\ 
 + 239764 n^8 - 18192 n^9 + 572 n^{10} + 4 n^{11})x_{13}^{8}\\
 -2 (-22336092 + 102242734 n - 207451554 n^2 + 246088757 n^3 - 189448366 n^4 + 99318803 n^5\\ 
 - 36163892 n^6 + 9130532 n^7 -1555886 n^8 + 166276 n^9 - 9008 n^{10} + 16 n^{11}\\ 
 + 16 n^{12})x_{13}^{7} + (-28792736 + 134723232 n - 268245792 n^2 + 301569274 n^3 \\ 
 -211579655 n^4 + 95750323 n^5 - 27336050 n^6 + 4215160 n^7 +516 n^8 - 146852 n^9\\ 
 + 30868 n^{10} - 2928 n^{11} + 112 n^{12})x_{13}^{6} -2 (65620080 - 227005040 n + 365106888 n^2 \\ 
 - 368314424 n^3 + 261802867 n^4 - 137589389 n^5 + 54260010 n^6 - 16033372 n^7\\ 
 + 3504636 n^8 - 550600 n^9 + 58748 n^{10} - 3800 n^{11} + 112 n^{12})x_{13}^{5}\\
 + (59501248 - 288258784 n + 627812916 n^2 - 804675780 n^3 + 674941445 n^4 - 391428642 n^5\\ 
 + 161902013 n^6 - 48456452 n^7 + 10471948 n^8 - 1597772 n^9 + 163380 n^{10} - 10040 n^{11}\\ 
 + 280 n^{12})x_{13}^{4} -4 (-24508224 + 70221360 n - 72172288 n^2 + 18690688 n^3\\
 + 27533464 n^4 - 33792759 n^5 + 19444979 n^6 - 7070399 n^7 + 1733171 n^8 \\
 - 287260 n^9 + 30980 n^{10} - 1968 n^{11} + 56 n^{12})x_{13}^{3} + 4 (-1 + n) (4378624 - 15912576 n \\
 + 22754384 n^2 - 15980992 n^3 + 4953656 n^4 + 505552 n^5 - 1046519 n^6 + 426038 n^7\\
 - 93847 n^8 + 12221 n^9 - 889 n^{10} + 28 n^{11})x_{13}^{2} -8 (-6 + n) (-4 + n) (-1 + n)^2 (2 + n) (31664 \\ 
 - 44256 n + 19472 n^2 - 1636 n^3 - 1423 n^4 + 535 n^5 - 76 n^6 + 4 n^7)x_{13}\\
 + 4 (-6 + n)^2 (-4 + n)^2 (-1 + n)^3 (2 + n)^2 (124 - 24 n - 5 n^2 + n^3)
\end{array}  
\end{equation*} }
For the case $x_{13} = 1$ system (\ref{eq11}) gives:
$$
f_{3} = x_{12}(x_{12}-x_{2}) = 0 \Longleftrightarrow x_{12} = x_{2}
$$
Thus we obtain the system
$$
x_{12} = x_{2},\ \ \ \ \ 3 x_2 + 2 x_{2}^2(2 - n) - x_{2}^3(1 - n) = 0
$$
The above system has three solutions:
$$
x_{2} = 0, \ \ x_{2} = \frac{-2 + n - \sqrt{7 - 7 n + n^2}}{n-1} \ \ x_{2} = \frac{-2 + n + \sqrt{7 - 7 n + n^2}}{n - 1}.
$$
The first solution is rejected, so we get the following
$$
x_{12} = x_{2} = \frac{-2 + n - \sqrt{7 - 7 n + n^2}}{n - 1}, \ \ x_{13} = x_{23} = 1
$$
and
$$
x_{12} = x_{2} = \frac{-2 + n + \sqrt{7 - 7 n + n^2}}{n - 1}, \ \ x_{13} = x_{23} = 1.
$$
These two solutions of system (\ref{eq11}) are Jensen's Einstein metrics on Stiefel manifolds pictured as follows:
$$
\begin{pmatrix}
 0 & \al &  1\\
\al & \al & 1\\
 1  & 1 &*
 \end{pmatrix}.
$$

Next we consider the case $x_{13}\neq 1$.  Then $h_{1}(x_{13}) = 0$ and we will   prove that the  equation $h_{1}(x_{13}) = 0$ has at least two positive roots.
We observe that 
$$
h_{1}(1) = 988524 - 7380396 n + 9224766 n^2 - 3877551 n^3 + 671409 n^4 - 
 40824 n^5
$$
is negative for $n\geq 7$ and
\begin{eqnarray*}
&& h_{1}(0) = -1142784 + 3459072 n - 3064576 n^2 - 222464 n^3 + 1556672 n^4\\
&& - 558592 n^5 - 114256 n^6 + 106352 n^7 - 19036 n^8 - 1072 n^9 + 
 776 n^{10} - 96 n^{11} + 4 n^{12}
\end{eqnarray*}
is positive for $n\geq 7$, hence we obtain one solution $x_{13} = \al_{13}$ between $0<\al_{13}<1$.  
In fact, it is possible to show that $1-4/n < \al_{13} <1-3/n$.
In the same way we observe that for $n\geq 7$ 
\begin{eqnarray*}
&& h_{1}(2) = -1077260544 + 2404260096 n - 1787496640 n^2 + 192056128 n^3 + 
 482885648 n^4 \\
 && - 324418304 n^5 + 97443600 n^6 - 15521168 n^7 + 
 1168004 n^8 \\
 && + 8544 n^9 - 8968 n^{10} + 528 n^{11} + 4 n^{12}
\end{eqnarray*}
is always positive,  hence we have a second solution $x_{13} = \beta_{13}$ between $1<\beta_{13}<2$.
In fact, it is possible to show that $1<\beta _{13}< 1+10/n^2$. 
\smallskip

Because for $x_{13} = 1$ we take  Jensen's metrics, we consider a Gr\"obner basis (take a lexicographic order $>$ with $z>x_{2}>x_{12}>x_{13}$) for the ideal $J$ generated by the polynomials
$
\{ f_1, \, f_2, \, f_3, $  $  \,z  \, x_2 \, x_{12} \, x_{13} \,( x_{13} - 1)  -1\}.
$
This basis contains the polynomial $h_{1}(x_{13})$ and the polynomials 
$$
6 (-6 + n) (-2 + n) (-1 + n)^3 (2 + n) (124 - 24 n - 5 n^2 + n^3) a(n)x_{2} - w_{2}(x_{13}, n),
$$
$$
12 (-4 + n) (-2 + n) (-1 + n)^2 a(n)x_{12} - w_{12}(x_{13}, n),
$$
where the first polynomial is obtained by the lexicographic order $z>x_{12}>x_{2}>x_{13}$ and the second by the lexicographic order $z>x_{2}>x_{12}>x_{13}$.
The polynomial $a(n)$ of $n$ is of degree $51$ and can be written as

{\small\small\small\begin{equation*} 
\begin{array}{l}
a(n) =1603489 (n-7)^{51}  +  331629 (n-7)^{50}  +
3330235 (n-7)^{49}  + 
 21654178 (n-7)^{48}  +\\
10256996 (n-7)^{47} 
 +  37749148 (n-7)^{46}
+  1124243 (n-7)^{45}  +  278618
(n-7)^{44} + \\
 586377 (n-7)^{43}+ 
10642404 (n-7)^{42}  +
16857968 (n-7)^{41}  +
23530271 (n-7)^{40}  +\\
29165722 (n-7)^{39}  +
32307839 (n-7)^{38}  +
32153102 (n-7)^{37}  +
28875189 (n-7)^{36}  +\\
23485774 (n-7)^{35}  +
17353528 (n-7)^{34}  +
11678023 (n-7)^{33}  +
71720699 (n-7)^{32}  +\\
40265082 (n-7)^{31}  +
20690586 (n-7)^{30}  +
97403757 (n-7)^{29}  +
42033064 (n-7)^{28}  +\\
16631548 (n-7)^{27}  +
60336057 (n-7)^{26} +
20061477 (n-7)^{25} +
61092907 (n-7)^{24} +\\
17022304 (n-7)^{23} +
43336722 (n-7)^{22} +
10063783 (n-7)^{21} +
21273147 (n-7)^{20} +\\
40831997 (n-7)^{19} +
70963420 (n-7)^{18} +
11131153 (n-7)^{17} +
15703148 (n-7)^{16} +\\
19849915 (n-7)^{15} +
22400611 (n-7)^{14} +
22495760 (n-7)^{13} +
20062962 (n-7)^{12} +\\
15890989 (n-7)^{11} +
11209194 (n-7)^{10} +
70758073 (n-7)^9 +
40094148 (n-7)^8 +\\
20263187 (n-7)^7 +
89201697 (n-7)^6 +
32751057 (n-7)^5 +
94463613 (n-7)^4 +\\
19935475 (n-7)^3 +
28717417 (n-7)^2 +
27565510 (n-7) +
16817600.
\end{array}  
\end{equation*} }

Hence, we see that for $n\geq 7$ the polynomial $a(n)$ is positive.  Thus for the positive values $x_{13} = \al_{13}, \beta_{13}$ found above, we obtain real values $x_{2} = \al_{2}, \beta_{2}$ and $x_{12} = \al_{12}, \beta_{12}$ as solutions of the system (\ref{eq11}).
We claim that $ \alpha_{2}, \beta_{2}, \alpha_{12}, \beta_{12}$ are positive.  We consider the ideal $J$ generated by 
$\{ f_1, \, f_2, \, f_3, $  $  \,z  \, x_2 \, x_{12} \, x_{13} \,( x_{13} - 1)  -1\}$ and now  take  a lexicographic order $>$  with $ z >   x_2 > x_{13} > x_{12}$ for a monomial ordering on $R$.  Then we see that a  Gr\"obner basis for the ideal $J$ contains the  polynomial $h_{2}(x_{12})$ 
{\small\small \begin{equation*} 
\begin{array}{l}
 h_{2}(x_{12}) = (-1 + n)^3 (-81 + 78 n - 22 n^2 + 2 n^3) (-1 + 22 n - 14 n^2 + 2 n^3)x_{12}^{10}\\
 -2 (-2 + n) (-1 + n)^2 (-5001 + 9686 n - 6508 n^2 + 1740 n^3 - 96 n^4 - 32 n^5 + 4 n^6)x_{12}^{9}\\
 + (-1 + n) (52393 - 76703 n - 6162 n^2 + 73700 n^3 - 56664 n^4 + 19352 n^5 - 3176 n^6 + 204 n^7)x_{12}^{8}\\
 -4 (-2 + n) (46818 - 149187 n + 191207 n^2 - 123037 n^3 + 41605 n^4 - 6890 n^5 \\
+ 394 n^6 + 10 n^7)x_{12}^{7} + (2108590 - 5378218 n + 5585179 n^2 - 2957751 n^3 + 809684 n^4\\ 
 - 92936 n^5 - 1860 n^6 + 912 n^7)x_{12}^{6} -2 (-2 + n) (829930 - 1018838 n \\ 
+ 348263 n^2 + 24142 n^3 - 31870 n^4 + 4152 n^5 + 16 n^6)x_{12}^{5} + (-150922 + 2293544 n \\ 
- 3272491 n^2 + 1786572 n^3 - 423924 n^4 + 33984 n^5 + 712 n^6)x_{12}^{4} -20 (-2 + n) (-104374 \\ 
+ 109331 n - 35890 n^2 + 2645 n^3 + 312 n^4)x_{12}^{3} + 25 (-76485 + 58035 n - 4052 n^2\\ 
 - 5488 n^3 + 1072 n^4)x_{12}^{2} -250 (-2 + n) (2051 - 1298 n + 224 n^2)x_{12} + 5625 (107 - 56 n + 8 n^2)
\end{array}  
\end{equation*} }
The polynomial  $h_{2}(x_{12})$ can be written as  
{\small\small\small\begin{equation*} 
\begin{array}{l}
h_{2}(x_{12}) = (2412504 + 5213484 (-7 + n) + 4966002 (-7 + n)^2 + 2728881 (-7 + n)^3\\
 + 950664 (-7 + n)^4 + 217336 (-7 + n)^5 + 32580 (-7 + n)^6 + 3088 (-7 + n)^7 \\
 + 168 (-7 + n)^8 + 4 (-7 + n)^9)x_{12}^{10}\\
 -(15481800 + 30522960 (-7 + n) + 26336250 (-7 + n)^2 + 13013966 (-7 + n)^3 + 4047916 (-7 + n)^4\\
  + 819792 (-7 + n)^5 + 107792 (-7 + n)^6 + 8840 (-7 + n)^7 + 408 (-7 + n)^8 + 8 (-7 + n)^9)x_{12}^{9}\\
  +(48245892 + 85904548 (-7 + n) + 66240641 (-7 + n)^2 + 28829438 (-7 + n)^3 + 7737264 (-7 + n)^4\\
   + 1310572 (-7 + n)^5 + 136796 (-7 + n)^6 + 8044 (-7 + n)^7 + 204 (-7 + n)^8)x_{12}^{8}\\
   -(97049440 + 154931088 (-7 + n) + 105623840 (-7 + n)^2 + 39883780 (-7 + n)^3\\
   + 9037872 (-7 + n)^4 + 1239340 (-7 + n)^5 + 97472 (-7 + n)^6 + 3736 (-7 + n)^7 + 40 (-7 + n)^8)x_{12}^{7}\\
   +(137946250 + 196719755 (-7 + n) + 117648788 (-7 + n)^2 + 38055081 (-7 + n)^3 \\
   + 7138384 (-7 + n)^4 + 767392 (-7 + n)^5 + 42828 (-7 + n)^6 + 912 (-7 + n)^7)x_{12}^{6}\\
   -(141888350 + 179760770 (-7 + n) + 93307270 (-7 + n)^2 + 25366350 (-7 + n)^3 \\
   + 3804144 (-7 + n)^4 + 298660 (-7 + n)^5 + 9808 (-7 + n)^6 + 32 (-7 + n)^7)x_{12}^{5}\\
   +(105439675 + 117258450 (-7 + n) + 51819665 (-7 + n)^2 + 11453180 (-7 + n)^3\\
    + 1288836 (-7 + n)^4 + 63888 (-7 + n)^5 + 712 (-7 + n)^6)x_{12}^{4}\\
    -(55868000 + 53548600 (-7 + n) + 19613300 (-7 + n)^2 + 3365760 (-7 + n)^3 \\
    + 258820 (-7 + n)^4 + 6240 (-7 + n)^5)x_{12}^{3}\\
    +(20567500 + 16633875 (-7 + n) + 4896700 (-7 + n)^2 + 613200 (-7 + n)^3 \\
    + 26800 (-7 + n)^4)x_{12}^{2} -(+4926250 + 3282750 (-7 + n) + 739500 (-7 + n)^2 + 56000 (-7 + n)^3)x_{12}\\
    +601875 + 315000 (-7 + n) + 45000 (-7 + n)^2
\end{array}  
\end{equation*} }

Then we have that, for $n\geq 7$ the coefficients of the polynomial $h_{2}(x_{12})$ is positive for even degree and negative for odd degree.  Thus if the equation $h_{2}(x_{12}) = 0$ has real solutions, then these is are all positive.
\par Next we take the lexicographic order $z>x_{12}>x_{13}>x_{2}$ for the monomial ordering $R$.  Then the Gr\"obner basis for the ideal $J$ contains the polynomial $h_{3}(x_{2})$ give by 
{\small\small\begin{equation*} 
\begin{array}{l}
h_{3}(x_{2}) = 81 (-1 + n)^3 (-7 + 3 n)^2 (124 - 24 n - 5 n^2 + n^3)x_{2}^{10} -54 (-2 + n) (-1 + n)^2 (-7\\
 + 3 n) (-12620 + 10066 n - 1796 n^2 - 63 n^3 + 20 n^4 + n^5)x_{2}^{9} \\
+ 9 (-1 + n) (5444348 - 12985508 n + 12566667 n^2 - 6115028 n^3 + 1510310 n^4 \\
- 165004 n^5 + 6142 n^6 - 961 n^7 + 153 n^8 + n^9)x_{2}^{8} -6 (-2 + n) (-2979568 \\
+ 8238150 n - 7797806 n^2 + 4055755 n^3 - 1567204 n^4 + 409230 n^5 \\
- 37660 n^6 - 4903 n^7 + 798 n^8 + 8 n^9)x_{2}^{7} +(664693040 - 1960669464 n + 2468056527 n^2\\
 - 1726283201 n^3 + 735938803 n^4 - 196313137 n^5 + 31516618 n^6 - 2509926 n^7 \\
+ 9408 n^8 + 8308 n^9 + 64 n^10)x_{2}^{6} -4 (-2 + n) (-59469494 + 203982575 n - 258060568 n^2\\
 + 163886493 n^3 - 58390388 n^4 + 11911514 n^5 - 1270014 n^6 + 44728 n^7 + 1634 n^8)x_{2}^{5}\\
 +(-1148666992 + 2483878280 n - 2018846305 n^2 + 682058094 n^3 + 1140539 n^4 - 74249792 n^5 \\
+ 24038016 n^6 - 3335112 n^7 + 180072 n^8)x_{2}^{4} -20 (-2 + n) (7338836 - 31781095 n + 35343394 n^2\\
 - 17774551 n^3 + 4654668 n^4 - 624056 n^5 + 34040 n^6)x_{2}^{3} + 100 (8425328 - 17876816 n \\
 + 15076868 n^2 - 6558218 n^3 + 1569983 n^4 - 197660 n^5 + 10272 n^6)x_{2}^{2}\\
 -1000 (-4 + n) (-2 + n) (-41149 + 30677 n - 7916 n^2 + 700 n^3)x_{2} \\
 + 22500 (-4 + n)^2 (107 - 56 n + 8 n^2)
\end{array}  
\end{equation*} }

The above polynomial can be rewritten as follows:
{\small\small\small\begin{equation*} 
\begin{array}{l}
h_{3}(x_{2}) = (185177664 + 353699136 (-7 + n) + 287249328 (-7 + n)^2 + 128537280 (-7 + n)^3 \\
+ 34614540 (-7 + n)^4 + 5746788 (-7 + n)^5 + 574533 (-7 + n)^6 + 31590 (-7 + n)^7 + 729 (-7 + n)^8)x_{2}^{10}\\
 -(1776660480 + 3383889696 (-7 + n) + 2750137920 (-7 + n)^2 + 1239055920 (-7 + n)^3\\
  + 338884128 (-7 + n)^4 + 57943890 (-7 + n)^5 + 6131268 (-7 + n)^6 + 381348 (-7 + n)^7 \\
  + 12420 (-7 + n)^8 + 162 (-7 + n)^9)x_{2}^{9}\\
  +(7529406192 + 14385411216 (-7 + n) + 11794825368 (-7 + n)^2 + 5412553560 (-7 + n)^3 \\
  + 1528684479 (-7 + n)^4 + 275559957 (-7 + n)^5 + 31789143 (-7 + n)^6 + 2286063 (-7 + n)^7\\
   + 96003 (-7 + n)^8 + 1998 (-7 + n)^9 + 9 (-7 + n)^10)x_{2}^{8}\\
  -(19003144320 + 36647336304 (-7 + n) + 30550728648 (-7 + n)^2 + 14409983892 (-7 + n)^3 \\
  + 4241037702 (-7 + n)^4 + 809574732 (-7 + n)^5 + 100606440 (-7 + n)^6 + 7901580 (-7 + n)^7 \\
  + 362442 (-7 + n)^8 + 8052 (-7 + n)^9 + 48 (-7 + n)^10)x_{2}^{7}\\
  +(31751845080 + 61819273140 (-7 + n) + 52334228666 (-7 + n)^2 + 25252421727 (-7 + n)^3 \\
  + 7654138784 (-7 + n)^4 + 1509829133 (-7 + n)^5 + 193077556 (-7 + n)^6 + 15306474 (-7 + n)^7\\
   + 673932 (-7 + n)^8 + 12788 (-7 + n)^9 + 64 (-7 + n)^10)x_{2}^{6}\\
  -(35544821200 + 69080883040 (-7 + n) + 58454603500 (-7 + n)^2 + 28210577088 (-7 + n)^3 \\
  + 8532473268 (-7 + n)^4 + 1666191528 (-7 + n)^5 + 207197760 (-7 + n)^6 + 15378664 (-7 + n)^7 \\
  + 577608 (-7 + n)^8 - 6536 (-7 + n)^9)x_{2}^{5}\\
  +(25602766500 + 48488384800 (-7 + n) + 39762417135 (-7 + n)^2 + 18448573770 (-7 + n)^3 \\
  + 5297021059 (-7 + n)^4 + 962339608 (-7 + n)^5 + 107676312 (-7 + n)^6 + 6748920 (-7 + n)^7 \\
  + 180072 (-7 + n)^8)x_{2}^{4}\\
  -(11214712000 + 19967074400 (-7 + n) + 15143699900 (-7 + n)^2 + 6348086000 (-7 + n)^3\\
   + 1588873060 (-7 + n)^4 + 237204560 (-7 + n)^5 + 19516480 (-7 + n)^6 + 680800 (-7 + n)^7)x_{2}^{3}\\
   +(2853346500 + 4609869000 (-7 + n) + 3090157200 (-7 + n)^2 + 1101382600 (-7 + n)^3 \\
   + 220180300 (-7 + n)^4 + 23376400 (-7 + n)^5 + 1027200 (-7 + n)^6)x_{2}^{2}\\
    -(387090000 + 547743000 (-7 + n) + 309590000 (-7 + n)^2 + 87525000 (-7 + n)^3 \\
    + 12384000 (-7 + n)^4 + 700000 (-7 + n)^5)x_{2}\\
    +21667500 + 25785000 (-7 + n) + 11587500 (-7 + n)^2 + 2340000 (-7 + n)^3 + 180000 (-7 + n)^4
\end{array}  
\end{equation*} }
Then we see that for $n\geq 7$ the coefficients of $h_{3}(x_{2})$ are positive for even degree and negative for odd degree.  Thus if the equation $h_{3}(x_{2}) = 0$ has real solutions then these are all positive.  Hence the numbers $\al_{12}, \al_{2}, \beta_{12}, \beta_{2}$ are positive.  In particular the positive solutions of (\ref{eq11}) are
$$
\{x_{2} = \al_{2}, x_{12} = \al_{12}, x_{13} = \al_{13}, x_{23} = 1\} \ \ \mbox{and} \ \
\{x_{2} = \be_{2}, x_{12} = \be_{12}, x_{13} = \be_{13}, x_{23} = 1\}
$$
and satisfy $\al_{13}, \be_{13}\neq 1$.  Thus, these solutions are different from Jensen's Einstein metrics, and can be pictured as
$$
\begin{pmatrix}
 0 & \be &  \gamma\\
\be & \al & 1\\
 \gamma  & 1 &*
 \end{pmatrix} \ \ \ (\al, \be, \gamma\ \ \mbox{are all different and}\ \ \gamma\neq 1).
$$

Similarly  we can show that the Stiefel manifold $V_5\mathbb{R} ^n = \SO(n)/\SO(n-5)$ with 
  $\Ad(\SO(2)\times\SO(3)\times\SO(n-5))$-invariant scalar products  of the form (\ref{metric1}) admits at least four invariant Einstein metrics. 
  Two of them are Jensen's Einstein metrics obtained before,  and  the other two are different from the Einstein metrics given by $\Ad( \SO(4)\times\SO(n-5))$-invariant scalar products. 
The computations are  quite long, so we present the proof very briefly.

We consider the system of equations 
  \begin{equation}\label{eq100}
   r_1 = r_2, \,  r_2 = r_{12}, \, r_{12} = r_{23},  \, r_{23}= r_{13}. 
  \end{equation} 
  We put $x_{23} = 1$

  { \small 
 \begin{equation}\label{}
\left. 
\begin{array}{lll} 
g_1 & = &  n {x_1} {x_{12}}^2
   {x_{2}}-n {x_{12}}^2
   {x_{13}}^2 {x_{2}}^2-5 {x_{1}}
   {x_{12}}^2 {x_{2}}+3 {x_{1}}
   {x_{13}}^2 {x_{2}}+5
   {x_{12}}^2 {x_{13}}^2
   {x_{2}}^2 \\
   & &
   +{x_{12}}^2
   -{x_{13}}^2-2
   {x_{13}}^2 {x_{2}}^2, \\
  g_2 & = &  n {x_{1}}
   {x_{12}}^2-n {x_{12}}^3
   {x_{13}}+n {x_{12}}
   {x_{13}}^3-2 n {x_{12}}
   {x_{13}}^2+n {x_{12}}
   {x_{13}}-5 {x_{1}}
   {x_{12}}^2 \\
   & &
+4 {x_{1}}
   {x_{13}}^2+5 {x_{12}}^3
   {x_{13}}-5 {x_{12}}
   {x_{13}}^3+4 {x_{12}}
   {x_{13}}^2-5 {x_{12}}
   {x_{13}}+2 {x_{13}}^2
   {x_{2}} \\
   g_3 & = & n {x_{1}} {x_{12}}^2-2
   n {x_{12}}^2 {x_{13}}-4
   {x_{1}} {x_{12}}^2+3 {x_{1}}
   {x_{13}}^2+3 {x_{12}}^3
   {x_{13}}+4 {x_{12}}^2
   {x_{13}} \\ & &-3 {x_{12}}
   {x_{13}}^3+3 {x_{12}}
   {x_{13}}, \\
   g_4 & = & 2 n {x_{12}}
   {x_{13}}^2-2 n {x_{12}}
   {x_{13}}+{x_{1}}
   {x_{12}}+{x_{12}}^2 {x_{13}}-2
   {x_{12}} {x_{13}}^2 {x_{2}}-4
   {x_{12}} {x_{13}}^2 \\ & &+4 {x_{12}}
   {x_{13}}-5 {x_{13}}^3+5
   {x_{13}}. 
\quad    
\end{array}  
\right\}
\end{equation}
}
 
  We consider a polynomial ring $R = {\mathbb Q}[z, x_1, x_2, x_{12}, x_{13}] $ and an ideal $I$ generated by $\{g_1, g_2, g_3, g_4, $ $ 
z  x_1  x_2  x_{12}  x_{13}-1\}$ to find  non zero solutions of the above equations. We take a lexicographic
order $>$ with $z > x_1 >  x_2 > x_{12} > x_{13}$ for a monomial ordering on $R$. Then, by the aid of computer,
we see that a Gr\"obner basis for the ideal $I$ contains the polynomial $(x_{13}-1) p_1(x_{13})$, where $ p_1(x_{13})$ is a polynomial of $x_{13}$ with  degree 22. 
By the same method as in the case $\Ad(\SO(4)\times\SO(n-5))$-invariant metrics, we can show that there are 
two positive
solutions of (\ref{eq100}) 
$$ \{ x_1 = \alpha_1, \, x_2 = \alpha_2, \, x_{12} = \alpha_{12} \, x_{13} = \alpha_{13}, x_{23} =1 \}, \{ x_1 = \beta_1, \, x_2 = \beta_2, \, x_{12} = \beta_{12} \, x_{13} = \beta_{13}, x_{23} =1 \}, 
$$
  with $ \alpha_{13} \neq 1, \beta_{13} \neq 1$ by considering the solutions  $ p_1(x_{13}) = 0 $.  
  In fact, we see that 
  \begin{eqnarray*}
  p_1(1) & = &  -51984 \left(5 n^3-44 n^2+130 n-75\right)^2 \left(56 n^5-921 n^4+5319 n^3-12654 n^2+10124 n-1356\right)\\
&=&  -51984 \left(5 (n-7)^3+61 (n-7)^2+249 (n-7)+394\right)^2 \times \\ & & (56 (n-7)^5+1039 (n-7)^4+6971
   (n-7)^3+20351 (n-7)^2+23529
   (n-7)+3754) < 0,  
 \end{eqnarray*}
 and also we see that 
 $p_1(1-5/n) > 0$   and $p_1(1+ 10/n^2) > 0$  for $n \geq 7$. Thus we see  $ 1- 5/n < \alpha_{13} < 1 $ and $ 1< \beta_{13} < 1+ 10/n^2$. 
 
These solutions  can be pictured as
$$
\begin{pmatrix}
 \al & \be &  \gamma\\
\be & \delta & 1\\
 \gamma  & 1 &*
 \end{pmatrix} \ \ \ (\al, \be, \gamma, \delta\ \ \mbox{are all different and}\ \ \gamma\neq 1).
$$ 
Note that in this case Jensen's metrics are pictured as
$$
\begin{pmatrix}
 \al & \al &  1\\
\al & \al & 1\\
 1  & 1 &*
 \end{pmatrix}.
$$

\end{proof}

\section{Generalized Wallach spaces}

A generalised Wallach space is a homogeneous space $M=G/K$ whose isotropy representation $\frak{m}$ decomposes into three $\operatorname{Ad}(K)$-invariant irreducible and pairwise orthogonal  submodules as
$\frak{m}=\frak{m}_1\oplus\frak{m}_2\oplus\frak{m}_3,
$
 which satisfy the relations
$[\frak{m}_i,\frak{m}_i]\subset \frak{k}$ $(i=1,2,3)$.  The original terminology of these spaces was 
three-locally-symmetric spaces (\cite{LNF}) since they generalize the defining property of classical symmetric spaces.

Some examples of generalised Wallach spaces are 
the Wallach spaces 
$\SU(3)/T_{\mbox{max}}$, $\Sp(3)/(\SU(2)\times \SU(2)\times \SU(2))$ and $F_4/\operatorname{Spin}(8)$,
the generalized flag manifolds 
$\SU(l+m+n)/\s(\U(l)\times \U(m)\times \U(n))$, $\SO(2l)/(\U(1)\times \U(l-1))$, $E_6/(\U(1)\times \U(1)\times \operatorname{Spin}(8))$, the homogeneous spaces
$\SO(l+m+n)/(\SO(l)\times \SO(m)\times \SO(n))$, $\Sp(l+m+n)/(\Sp(l)\times \Sp(m)\times \Sp(n))$, and (as special case)
the Stiefel manifolds  $\SO(n+2)/\SO(2)$.
As a consequence of their definition if follows that $[\frak{m}_i,\frak{m}_j]\subset \frak{m}_k$, for $i, j, k$ distinct. 
Despite their simple description, a complete  classification of generalized Wallach spaces was given only  recently by Yu.G. Nikonorov in \cite{Ni2} and Z. Chen, Y. Kang and K. Liang in \cite{ChKaLi}.

Invariant Einstein metrics on generalized Wallach spaces were been originally studied in \cite{Ni1}.
In that paper it is proved that every generalized Wallach space admits at least one invariant Einstein metric, and in \cite{LNF}  that it admits at most four invariant Einstein metrics (up to a homothety).
A good survey about them can be found in \cite{ChKaLi} and \cite{ChNi}.

The most subtle example is the homogeneous space $\SO(l+m+n)/(\SO(l)\times \SO(m)\times \SO(n))$,
and in \cite{ChNi} there is serious progress towards the classification of Einstein metrics in this space.
It turns out that the number of invariant Einstein metrics on  $\SO(l+m+n)/(\SO(l)\times \SO(m)\times\SO(n))$
could be estimated by using  special properties of the normalized Ricci flow on generalized Wallach spaces (cf. \cite{AbiArNiSi}, \cite{ChNi}).
The main results in \cite{ChNi} are the following:

\begin{theorem}
Assume that $l\geq m \geq n \geq 1$ and $m\geq 2$. Then
the number of invariant Einstein metrics on the space $G/K=\SO(l+m+n)/(\SO(l)\times \SO(m)\times \SO(n))$ is
four for $n> \sqrt{2l+2m-4}$, and two for $n<\sqrt{l+m}$ (up to a homothety).
\end{theorem}

\begin{theorem}
Let $q =2, 3$ or $4$.  Then  there are infinitely many homogeneous spaces $\SO(l + m + n)/(\SO(l) \times \SO(m) \times \SO(n))$ that admit exactly $q$ invariant Einstein
metrics up to a homothety.
\end{theorem}

  \section{Some open problems}
  There is no doubt that there has been a lot of progress towards the classification of invariant Einstein metrics on generalized flag manifolds.  It seems that Ziller's finiteness conjecture is in fact true in this case, that is the number of Einstein metrics (up to isometry) is finite.  However, as the number of isotropy summands increases then  determining the total number of Einstein metrics is getting a difficult task, especially for flag manifolds determined by classical Lie groups.  The Gr\"obner bases techniques used can contribute a lot, but  seem to be unable to shed more light in the set of invariant Einstein metrics.  It seems that new tools and techniques have to be used.  The deep works \cite{Bom} by C. B\"ohm and \cite{BWZ} by C. B\"ohm, M. Wang and W. Ziller constitute an alternative point of view.  Also, the recent works by M.M. Graev \cite{Gr1},  \cite{Gr2} and \cite{Gr3} are important  contributions towards the understanding of the number of complex Einstein metrics on flag manifolds.  
  
  Concerning invariant Einstein metrics on Stiefel manifolds (or other homogeneous spaces where the isotropy representation contains equivalent summnands) the picture is more foggy. 
The set of invariant  metrics is quite vast in this case, so searching for invariant Einstein metrics is not an easy matter.  One has to search in certain subsets of the set of invariant metrics (``symmetry assumption") so that computations of the Ricci tensor get slightly simpler, and then hope to find Einstein metrics in these subsets.  
  Even though there is considerable success by such approach, there is no particular evidence on whether the finiteness conjecture is true or not.  There might even exist  one-parameter families of invariant Einstein metrics.
  
  Finally some open problems related to the generalized Wallach space    $\SO(l + m + n)/(\SO(l) \times \SO(m) \times \SO(n))$ are listed in \cite{ChNi}.

  \medskip

\noindent
{\bf Acknowledgements.} 
The work was supported by Grant $\# E.037$ from the Research Committee of the University of Patras
(Programme K. Karatheodori).
I would like to take the opportunity to express my gratitude to Professor M. Guest for directing me in various areas of geometry including homogeneous geometry.  Through the years I had (and still have) beneficial  collaborations with   D.V. Alekseevsky, Yu. G. Nikonorov, Y. Sakane, I. Chrysikos, M. Statha and N.P. Souris.  It is a pleasure to still have fruitful discussions with them.

\end{document}